\documentclass[a4paper]{amsart}
\input{commandoos.sty}

\begin{document}
\author{A. Eggink}
\title{Hilbert's Tenth Problem for some Noncommutative Rings} 
\begin{abstract}
	We consider Hilbert's tenth problem for two families of noncommutative rings. Let $K$ be a field of characteristic $p$. We start by showing that Hilbert's tenth problem has a negative answer over the twisted polynomial ring $K\{\t\}$ and its left division ring of fractions $K(\t)$. We prove that the recursively enumerable sets and Diophantine sets of $\F_{q^n}\{\t\}$ coincide. We reduce Hilbert's tenth problem over $\F_{q^n}\{\!\{\t\}\!\}$ and $\F_{q^n}(\!(\t)\!)$, the twisted version of the power series and Laurent series, to the commutative case. Finally, we show that the different models of $\F_q[T]$ in $K\{\t\}$ we created are all equivalent in some sense which we will define. We then move on to the second family of rings, coming from differential polynomials. We show that Hilbert's tenth problem over $K[\dd]$ has a negative answer. We prove that Hilbert's tenth problem over the left division ring of fractions $K(\dd_1,\ldots,\dd_k)$ can be reduced to Hilbert's tenth problem over $C(t_1,\ldots,t_k)$ where $C$ is the field of constants of $K$. This gives a negative answer for $k \geq 2$ if the field of constants is $\C$ and for $k\geq 1$ if it is $\R$. 
\end{abstract}

\maketitle 

\section{Introduction}

\subsection{Hilbert's tenth problem} \phantom{=}\\
In 1900 David Hilbert posed a list of problems. The tenth one asks for an algorithm that given a Diophantine equation with integer coefficients decides whether or not that Diophantine equation has a solution in the integers. It turns out that such an algorithm does not exist, which was proven in 1970 by Yuri Matiyasevich, using the work of Julia Robinson, Martin Davis and Hilary Putnam. \\
Hilbert's tenth problem can easily be generalized, by replacing the integers with another ring in which the solutions must lie and possibly replacing the ring of coefficients. Another way of phrasing this is asking whether the positive existential theory over a certain ring is decidable. This generalized version has been considered before. The most important open question is that of Hilbert's tenth problem over $\Q$. Hilbert's tenth problem has a negative answer over $\F_q(T)$ when $T$ is added to the language proven by Pheidas \cite{PheidasFpT} and Videla \cite{H10PF2(T)}, over $\C(t_1,t_2)$ when $t_1$ and $t_2$ are added to the language by Kim and Roush \cite{KRC(t1t2)} and over $R[T]$ when $T$ is added to the language for a commutative domain $R$ by Denef \cite{DenefRTchar0}, \cite{DenefRTcharp} and there are many more results.

\subsection{New results} \phantom{=}\\
The main results of this paper are about two noncommutative families of rings. The first is that of the twisted polynomial rings and variants, which we define below. 

\begin{defi}
	Take $K$ a field of characteristic $p$ and let $q$ be a power of $p$. We define $K\{\t\}$ to be the \textbf{twisted polynomial ring}, with the usual addition but the multiplication defined by $\t a = a^q \t$ for all $a \in K$. \\
	We define $K(\t)$ to be the left division ring of fractions of $K\{\t\}$, which exists since $K\{\t\}$ satisfies the left Ore condition, see Exercise 4.1.9 of \cite{DM}. \\
	The rings $K\pownon{\t}$ and $K\lau{\t}$ are the twisted version of the formal power series and the formal Laurent series, respectively. 
\end{defi} 

\begin{thm} \phantom{=}
	\begin{enumerate}
		\item Hilbert's tenth problem has a negative answer over $K\{\t\}$ with ring of coefficients $\F_p\{\t\}$.
		\item Hilbert's tenth problem has a negative answer over $K(\t)$ with ring of coefficients $\F_p(\t)$.
		\item If Hilbert's tenth problem has a negative answer over $\F_q\pow{T}$ or $\F_q\lau{T}$ with ring of coefficients $\F_p\pow{T}$ or $\F_p\lau{T}$, then it has so over $K\pownon{\t}$ or $K\lau{\t}$ with ring of coefficients $\F_p\pownon{\t}$ and $\F_p\lau{\t}$, respectively.
		\item If Hilbert's tenth problem has a positive answer over $\F_p\pow{T}$ or $\F_p\lau{T}$ with ring of coefficients $\F_p\pow{T}$ or $\F_p\lau{T}$, then it has so over $\F_{q^n}\pownon{\t}$ or $\F_{q^n}\lau{\t}$ with ring of coefficients $\F_q\pownon{\t}$ or $\F_q\lau{\t}$, respectively. 
	\end{enumerate}
\end{thm}
We note that under the condition that we can resolve singularities in characteristic $p$, it is proven in \cite{Fq[[T]]DenefSchoutens} that Hilbert's tenth problem over the rings $\F_q\pow{T}$ and $\F_q\lau{T}$ has a positive answer. If the resolution of singularities in characteristic $p$ is true, the third case in the theorem does not occur.\\
We have one stronger result, namely that the Diophantine sets (for the definition see \th\ref{toolsDefDiophSet}) and the recursively enumerable sets (for the definition see Appendix A of \cite{Shlapentokh}) of the ring $\F_{q^n}\{\t\}$ are the same if we take the coefficient ring to be the ring itself. \\
The last result in this setting is that the different Diophantine models of $\F_q[T]$ in $K\{\t\}$ we create are equivalent. Intuitively a Diophantine model of $R_0$ in $R_1$ is a subset of a ring $R_1$ that contains the Diophantine structure of $R_0$. The intuition for two models being equivalent is that  we can map one to the other in a ``Diophantine" way. We will properly define these notions in Paragraph \ref{parDMModels}.  \\

The second family of rings we consider are the rings of differential polynomials and their left division rings of fractions. This is covered in Section \ref{chapterDiff}. \\
We first recall what differential rings are, see also Paragraph 1.1 of \cite{Kolchin} :
\begin{defi}
	A \textbf{derivation} on a commutative ring $R$ is a function $\dd:R \ra R$ such that for all $r,s \in R$ it holds that $\dd(r+s) = \dd(r) + \dd(s)$ and $\dd(rs) = r \dd(s) + s\dd(r)$.\\
	A \textbf{differential ring} is a commutative ring $R$ with derivations $\dd_1,\ldots,\dd_k$ such that $\dd_i(\dd_j(r)) = \dd_j(\dd_i(r))$ for all $i,j \in \{1,\ldots,k\}$ and $r \in R$. \\ 
	A \textbf{differential field} is a differential ring $K$ that is also a field.\\
	The set $\{f \in R\mid \dd_1(f)=\ldots=\dd_k(f) = 0\}$ is the \textbf{ring or field of constants} of the differential ring or field $R$ that has derivations $\dd_1,\ldots,\dd_k$.
\end{defi} 
Before we define differential polynomials I want to remark two things. First, in what follows I use $\frac{d}{dx}$ and $\frac{d}{dx_i}$ for arbitrary derivations if $x$ and $x_i$ are not defined. Second, there are two versions of the ring of differential polynomials, a commutative and a noncommutative version. I use the noncommutative version, which is sometimes also called a Weyl algebra, see also \cite{WikiWeyl} and Example 1.9 of \cite{Lam}. The definition in Paragraph 1.6 of \cite{Kolchin} is different. 
\begin{defi} 
	Let $K$ be a differential ring or field with derivations $\frac{d}{dx_1},\ldots,\frac{d}{dx_k}$. The ring of \textbf{differential polynomials} is the set $K[\dd_1,\ldots,\dd_k]$ where the addition is the usual addition, but the multiplication is defined by $\dd_i\dd_j = \dd_j\dd_i$  and $\dd_i a = a \dd_i + \frac{d}{dx_i}(a)$ for all $i,j \in \{1,\ldots,k\}$ and all $a \in K$. \\
	The ring $K(\dd_1,\ldots,\dd_k)$ is the left division ring of fractions of $K[\dd_1,\ldots,\dd_k]$, which exists since $K[\dd_1,\ldots,\dd_k]$ satisfies the left Ore condition, see \th\ref{diffOreCondition}. 
\end{defi}

\begin{thm}
	Let $(K,\frac{d}{dx})$ be a differential field. Then Hilbert's tenth problem over $K[\dd]$ with ring of coefficients $\Q[\dd]$ or $\F_p[\dd]$ has a negative answer.
\end{thm}

\begin{thm}
	Let $K$ be a differential field with field of constants $C$. Hilbert's tenth problem has a negative answer over $K(\dd_1,\ldots,\dd_k)$ with ring of coefficients $\Q(\dd_1,\ldots,\dd_k)$ if $C = \R$ or if $C = \C$ and $k \geq 2$.  
\end{thm}

\subsection{Methods}\phantom{=}\\
The method we use to prove this results is modeling a ring for which we know that Hilbert's tenth problem has a negative answer inside the ring we are interested in. I formalized this common practice in \cite{ArtikelDioph} with the notion of an (effective) Diophantine map. We will introduce the main results we need in Section \ref{chapterTools}. \\
The map we use for the twisted polynomial ring is the Carlitz module:
\begin{defi}
	Set $A = \F_q[T]$. Let $\y: A \ra K$ be an $\F_q$-algebra homomorphism, which makes $K$ into an $A$-field. Then the \textbf{Carlitz module} is the $\F_q$-algebra homomorphism defined by $\phi(T) = \phi_T = \y(T) + \t$. See Definition 3.2.5 in \cite{DM}.   
\end{defi}
The maps we use for the variants of the twisted polynomials are all extensions of the Carlitz module. Note that the Carlitz module is a special case of a \textbf{Drinfeld module}, which is an $\F_q$-algebra homomorphism $\phi: \F_q[T] \ra K\{\t\}$ defined by $T \ra \y(T) + g_1\t + \ldots + g_r\t^r$, see Definition 3.2.2 of \cite{DM}.\\
For the case of the differential polynomials we use the $C$-algebra homomorphisms $\phi: C[t] \ra K[\dd]$ defined by $\phi(t) = \dd$ and $\psi: C(t_1,\ldots,t_k) \ra K(\dd_1,\ldots,\dd_k)$ defined by $\psi(t_i) = \dd_i$. The first map is a special case of a Krichever module:
\begin{defi}
	Let $(K,\frac{d}{dx})$ be a differential field with field of constants $\C$. A \textbf{Krichever module} is a $\C$-algebra homomorphism $\phi: \C[t] \ra K[\dd]$ such that $\phi(\C[t]) \not \subset K$. 
\end{defi} 
Note that this is very similar to a Drinfeld module, since both are algebra homomorphisms whose image should not be contained in the base field. For a Drinfeld module, we have the extra requirement that it respects the $\y$-structure of $K$, which we do not have in the setting of differential polynomials.

\section{Acknowledgements} 
I thank my supervisor, Gunther Cornelissen, for suggesting the topic, his supervision during the project and his help in writing this paper.

\section{Conventions}
The rings have by assumption a $1$, but they do not have to be commutative. Furthermore, I assume that each language contains a constant $0$ and a relation ``$=$", where often ``$=$" is left out from the notation.\\

\section{Notation} 
\begin{itemize}
	\item $\N = \Z_{\geq 0}$;
	\item $\L_R = \{0,1,+,\cdot\}$ is the language of rings;
	\item $p_n$ is the $n$-th prime number;
	\item $\rho(i,n)$ gives the exponent of $p_i$ in the prime factorization of $n$.
	\item $(\mu x)(R(x))$ is the minimizing operator, which gives the least $x$ such that the relation $R(x)$ is true;
	\item $S_c(\L)$ is the set of constants of a language $\L$;
	\item $S_c'(\L) = \{x\in R\mid \{x\}\text{ is a Diophantine set}\}$, for a language $\L$ with an interpretation in the set $R$; 
	\item $\L_t= \{0,1,t,+,\cdot\}$ and $\L_T = \{0,1,T,+,\cdot\}$;
	\item $\L_{t_1,t_2} = \{0,1,t_1,t_2,+,\cdot\}$;
	\item $\L_{\t,\y} = \{0,1,\t,\y(T),+,\cdot\}$ in Section \ref{chapterDM};
	\item $\L_{\dd} = \{0,1,\dd,+,\cdot\}$;
	\item $\L^{(k)}$ is the language on $R^k$ that creates the same Diophantine sets as $(R,\L)$, see \th\ref{DMdefcartlanguage};
	\item $A = \F_q[T]$ and $F = \F_q(T)$, only in Section \ref{chapterDM};
	\item $\y:A \ra K$ is the $A$-field structure of $K$, only in Section \ref{chapterDM};
	\item $\a$ is a generator of $\F_{q^n}$ over $\F_q$, only in Paragraphs \ref{parDMMRDP},\ref{parDMpower} and \ref{parDMlaurent};
	\item $\b$ is a generator of $\F_q$ over $\F_p$, only in Paragraphs \ref{parDMMRDP},\ref{parDMpower} and \ref{parDMlaurent}.
	\item $K_l = \{c \in K\mid \frac{d}{dx_i} c = 0\}$, only in Section \ref{chapterDiff};
	\item $R_{l} = K_l(\dd_1,\ldots,\dd_{l-1},\dd_{l+1},\ldots,\dd_k)[\dd_{l}]$, only in Section \ref{chapterDiff};
	\item $\dd^{\vec{i}} = \dd_1^{i_1}\dd_2^{i_2}\cdots\dd_k^{i_k}$;
	\item $\frac{d}{dx_j} \sum_{\vec{i}} a_{\vec{i}}\dd^{\vec{i}} = \sum_{\vec{i}} \frac{d}{dx_j}(a_{\vec{i}}) \dd^{\vec{i}}$ for a derivation $\frac{d}{dx_j}$ and an element $\sum_{\vec{i}} a_{\vec{i}} \dd^{\vec{i}} \in K[\dd_1,\ldots,\dd_k]$.
\end{itemize}

\section{Tools} \label{chapterTools} 
In this section we gather the most important tools of \cite{ArtikelDioph} that we need. We also use some theory of recursive functions, recursive sets, recursively enumerable sets and recursive rings. The theory of this can be found in Appendix A of \cite{Shlapentokh}. \\
We start with a few definitions:

\begin{defi} \th\label{toolsDefDiophSet} 
	Let $(R,\L)$ be a set and a language (with an interpretation in $R$). A \textbf{term} of $\L$ consists of constants of $\L$ and variables, combined with functions of $\L$. A \textbf{basic formula} is of the form $(t_1,\ldots,t_k)\in R'$ where $R' \in \L$ is an $k$-ary relation (which can be equality) and the $t_i$ are terms of $\L$.\\
	A set $S \subset R^k$ is Diophantine if there exist basic formulas $f_1,\ldots,f_r$ such that
	$$S = \{(x_1,\ldots,x_k) \in R^k\mid \exists y_1,\ldots,y_l \in R\ (f_1\wedge f_2\wedge \ldots \wedge f_r)\}.$$   
\end{defi}
If we use the language of rings $\L_R = \{0,1,+,\cdot\}$ we get the more common definition of a Diophantine set back.

\begin{notation}
	Let $(R,\L)$ be a set with a language. We denote with $S_c(\L)$ the set of constants of $\L$ and we let $S_c'(\L) = \{x \in R\mid \{x\} \text{ is a Diophantine set}\}$.
\end{notation}

\begin{defi}
	Let $(R_1,\L_1)$ and $(R_2,\L_2)$ be sets with a language. A map $d: R_1 \ra R_2$ is called a \textbf{Diophantine map} if for all Diophantine sets $S \subset R_1^k$ the set $d(S) \subset R_2^k$ is Diophantine.
\end{defi}
In this article, we will use \textbf{effective Diophantine maps}. We will not give a definition here. Intuitively it is a Diophantine map with the extra property that we can construct the Diophantine equations that describe $d(S)$ from the ones that describe $S$. \\

There are some technicalities when defining Hilbert's tenth problem since algorithms cannot handle an uncountable amount of options. This leads to the following definition:
\begin{defi} \th\label{H10PexistenceH10P}
	Take $(R,\L)$ to be a set with a language. We say that \textbf{Hilbert's tenth problem exists over $(R,\L)$} if there exists a recursive ring $R_0$ such that $S_c(\L) \subset R_0 \subset S_c'(\L)$.
\end{defi}

With this, we have gathered the definitions to state the two theorems we will use:
\begin{thm}\th\label{H10PR2deci+deff=>R1deci}
	Let $(R_1,\L_1)$ and $(R_2,\L_2)$ be sets with a countable language. Suppose that Hilbert's tenth problem exists over $(R_1,\L_1)$ and $(R_2,\L_2)$. If Hilbert's tenth problem over $(R_2,\L_2)$ has a positive answer and $d:R_1 \ra R_2$ is an effective Diophantine map, then Hilbert's tenth problem over $(R_1,\L_1)$ has a positive answer.
\end{thm}

\begin{thm} \th\label{H10PeffectiveDioph} 
	Take $(R_1,\L_1)$ and $(R_2,\L_2)$ sets with a finite language. Suppose that Hilbert's tenth problem exists over $(R_1,\L_1)$ and $(R_2,\L_2)$. Let $d: (R_1,\L_1) \ra (R_2,\L_2)$ be an injective map. 
	Suppose that for all the following sets:
	\begin{itemize}
		\item For all constants $c \in \L_1$, the set $\{d(c)\}$;
		\item For all functions $f: R_1^k \ra R_1$ of $\L_1$ the set $\{(d(x_1),\ldots,d(x_k),d(f(x_1,\ldots,x_k))) \mid \vec{x} \in R_1^k\}$;
		\item For all relations $S \subset R_1^k$ of $\L_1$ the set $d(S)$;
		\item The set $d(R_1)$,
	\end{itemize}
	an explicit description of the form 
	\begin{align} 
		\{\vec{x} \in R_2^k \mid \ex \vec{y} \in R_2^l,\ \va 1\leq i \leq n,\ (t_{i1}(\vec{x},\vec{y}),\ldots,t_{in_i}(\vec{x},\vec{y})) \in S_i\} \label{eqExplFormDioph}
	\end{align} 
	with $t_{ij}$ terms of $\L_2$ and $S_i$ relations of $\L_2$ is known.\\
	Then the map $d$ is an effective Diophantine map.
\end{thm} 

\section{Hilbert's Tenth Problem over the Twisted Polynomial Ring and Variants}  \label{chapterDM}
In this section, we will first prove that Hilbert's tenth problem over the twisted polynomial ring $K\{\t\}$ has a negative answer. We extend the negative answer to Hilbert's tenth problem over $K\{\t\}$ to $K(\t)$, the left division ring of fractions of $K\{\t\}$. Then we specialize to $K= \F_{q^n}$ and we prove that MRDP holds in $\F_{q^n}\{\t\}$. Next, we take a look at some reductions for Hilbert's tenth problem over $\F_{q^n}\pownon{\t}$ and $\F_{q^n}\lau{\t}$. Finally, we show that we can create different models of $\F_q[T]$ in $K\{\t\}$ and we prove that they are all equivalent.\\
In this section we make use of Drinfeld modules, the theory of this can be found in \cite{DM}. We will mostly use the results and definitions of Chapter 3 in \cite{DM}.  

\subsection{A Negative Answer over $K\{\t\}$} \label{parDMKt} \phantom{=}\\
We will show in this paragraph that the Carlitz module $\phi: \F_q[T] \ra K\{\t\}$ is an effective Diophantine map and we will conclude that Hilbert's tenth problem over $K\{\t\}$ has a negative answer.

\begin{defi}[Definition 3.3.1 of \cite{DM}]
	The endomorphisms of $\phi$ are
	$$\End_K(\phi) = \{u \in K\{\t\}\mid u\phi_T = \phi_T u\}.$$
\end{defi} 

\begin{lemma} \th\label{DMphiAeq} 
	We have $\End_K(\phi)=\phi(A)$
\end{lemma} 

\begin{proof}
	By Theorem 3.4.1 of \cite{DM} we get that $\End_K(\phi)$ is a free $A$-module of rank at most $1$, since the Carlitz module has rank $1$. We have that $1 \in \End_K(\phi)$, so the rank is $1$. Since the $A$-module structure of $\End_K(\phi)$ is given by multiplication by $\phi_a$, this implies that $\phi(A)= \End_K(\phi)$.
\end{proof} 

\begin{notation}
	We denote with $\L_{\t,\y} = \{0,1,\y(T),\t,+,\cdot\}$ the language of rings with extra constants $\t$ and $\y(T)$. 
\end{notation} 

\begin{remark} \th\label{DMdefy0} 
	Let $\y: A \ra K$ be an arbitrary $A$-field structure. We can define a new $A$-field structure $\y_0$ by $\y_0(f(T)) = \y(f(0))$, so we first reduce modulo $T$, include the result in $A$ and use $\y$ to map to $K$. Note that $\y_0(T) = \y(0) = 0$, so by choosing an appropriate $A$-field structure, we can remove the constant $\y(T)$.   
\end{remark}

\begin{thm} \th\label{DMphiAdioph} 
	Let $\L_T$ be the language on $A$ and let $\L_{\t,\y}$ be the language on $K\{\t\}$. Then the map $\phi:A \ra K\{\t\}$ is an effective Diophantine map.
\end{thm} 

\begin{proof}
	The map $\phi$ is injective. This means we are going to use \th\ref{H10PeffectiveDioph}. 
	\begin{itemize}
		\item By \th\ref{DMphiAeq} 
		$$\phi(A) = \End_K(\phi) = \{u \in K\{\t\}\mid u\phi_T = \phi_T u \}  = \{u\in K\{\t\} \mid u(\y(T) + \t) = (\y(T) + \t)u\},$$
		which is a Diophantine set.
		\item The elements $\phi(0) = 0$, $\phi(1) = 1$ and $\phi(T) = \y(T) + \t$ are all Diophantine expressions in the constants of $\L_{\t,\y}$.
		\item Let $* \in \{+,\cdot\}$. Since $\phi$ is a ring homomorphism, we get
		\begin{align*} 
			\phi(\{(a,b,a*b)\mid a,b \in A\}) = \{\{\phi(a),\phi(b),\phi(a)*\phi(b)\mid a,b \in A\} \\
			=\{(x,y,z)\in K\{\t\} \mid x,y \in \phi(A),\ z = x*y\}.
		\end{align*} 
		Since $\phi(A)$ is a Diophantine set, we get that the set $\phi(\{(a,b,a*b)\mid a,b \in A\})$ is Diophantine.
	\end{itemize}
	We have checked all conditions of \th\ref{H10PeffectiveDioph}, so $\phi$ is an effective Diophantine map.
\end{proof} 

\begin{lemma} \th\label{DMrecreprKtau} 
	Let $L$ be a field of characteristic $p$ with a recursive representation. Then $L\{\t\}$ has a recursive representation. Furthermore there exist a recursive function $C: L\{\t\} \ti \N \ra K$ such that $C(\sum_{i=0}^n a_i\t^i,k) = a_k$ and a recursive function $\deg_\t:L\{\t\} \ra \N$ that gives the degree in $\t$. 
\end{lemma}

\begin{proof}
	Let $j$ be a recursive representation of $L$ such that $j(0) = 0$. Such a representation exists since we can compose a recursive representation with the recursive function $f: \N \ra \N$ that swaps the values $j(0)$ and $0$. Define $J(\sum_{i=0}^na_i\t^i) = \prod_{i=0}^n p_{i+1}^{j(a_i)}$, where $p_i$ is the $i$-th prime number. $J$ is well-defined, since $p_{i+1}^{j(a_i)} = 1$ if $a_i=0$. Furthermore $J$ is injective, since $j$ is injective.\\
	Let $l(n)$ be the largest prime dividing $n$ and let $\rho(i,n)$ be the exponent of $p_i$ in the prime factorization of $n$. Set $\X_S$ to be the characteristic function of $S$. Then 
	$$J(L\{\t\}) = \{ n \in \N\mid 1 = \prod_{i=0}^{l(n)} \X_{j(L)}(\rho(i,n)) \},$$ which is a recursive set since $j(L)$ is a recursive set. Let $P_+,P_-$ and $P_\ti$ be the translations of $+$, $-$ and $\cdot$ under $j$, i.e.\ $P_*(n,m) = j(j^{-1}(n)*j^{-1}m)$. Then the translation of $+$ under $J$ is
	$$T_+(n,m) = \prod_{i=0}^{\max(l(n),l(m))} p_{i+1}^{P_+(\rho(i,n),\rho(i,m))},$$
	which is a recursive function. In the same way the translation $T_-$ of $-$ is a recursive function.\\
	We have that $$\left(\sum_{i=0}^na_i\t^i\right)\left(\sum_{i=0}^mb_i\t^i\right) = \sum_{k=0}^{m+n} \left(\sum_{i=0}^k a_ib_{k-i}^{q^i}\right)\t^k.$$
	This gives for the translation of multiplication:
	$$T_{\ti}(n,m) = \prod_{k=0}^{l(n)+l(m)} p_{k+1}^{(\sum')_{i=0}^k (P_\ti(\rho(i,n),(\prod')_{k'=0}^{q^i} \rho(k-i,m)))}.$$
	Here the $\sum'$ and $\prod'$ need to be understood using $P_+$ and $P_\ti$ instead of $+$ and $\cdot$. This gives that $T_\ti$ is a recursive function.\\
	We checked all conditions, guaranteeing that $J$ is a recursive representation of $L\{\t\}$.\\
	We have that $j(C(J^{-1}(n),k)) = \rho(k,n)$ and $\deg_\t(J^{-1}(n)) = l(n)$, which are both recursive functions. 
\end{proof}

\begin{thm} \th\label{DMKtauUndecidable} 
	Hilbert's tenth problem over $(K\{\t\},\L_{\t,\y})$ exists and has a negative answer.
\end{thm}

\begin{proof}
	By \cite{DenefRTcharp} Hilbert's tenth problem over $A = \F_q[T]$ with the language $\L_T$ has a negative answer.
	The ring $\F_p(\y(T))$ is a recursive ring, so with \th\ref{DMrecreprKtau} we get that $\F_p(\y(T))\{\t\}$ has a recursive representation. Since we have that $S_c(\L_{\t,\y}) \subset \F_p(\y(T))\{\t\} \subset S_c'(\L_{\t,\y})$, this implies that Hilbert's tenth problem exists over $(K\{\t\},\L_{\t,\y})$. \\
	With \th\ref{DMphiAdioph} we get that $\phi: A \ra K\{\t\}$ is an effective Diophantine map. Then \ref{H10PR2deci+deff=>R1deci} implies that Hilbert's tenth problem over $(K\{\t\},\L_{\t,\y})$ has a negative answer.
\end{proof} 

\subsection{A Negative Answer for $K(\t)$} \phantom{=}\\
In the previous paragraph, we have seen that $A$ has a Diophantine model in $K\{\t\}$. The next logical question is if we can embed the fraction field $F$ into some kind of fraction field of $K\{\t\}$. Since $K\{\t\}$ satisfies the left Ore condition, see exercises 9 of paragraph 4.1 of \cite{DM}, we get a left division ring of fractions $K(\t)$. In this paragraph, we will show that $\phi: F \ra K(\t)$, the extension of $\phi:A \ra K\{\t\}$, is an effective Diophantine map. From this, we will conclude that Hilbert's tenth problem over $K(\t)$ has a negative answer. 

\begin{prop} \th\label{DMcoefpoly} 
	Let $f_1,f_2,g \in K\{\t\}$ and suppose that $g$ is nonzero. Write $m = \deg_\t(g)$ and $n = \max(\deg_\t(f_1),\deg_\t(f_2))$. We can write $f_1 = \sum_{i=1}^{n} a_i\t^i$, $f_2 = \sum_{i=1}^{n} b_i\t^i$, $f_1g = \sum_{i=0}^{n+m} c_i\t^i$ and $f_2g = \sum_{i=0}^{n+m}d_i\t^i$. If $c_i=d_i$ for all $i\geq N$, then $a_i=b_i$ for all $i \geq N-m$. \\
	The same holds if $gf_1$ and $gf_2$ have the same coefficients at $\t^i$ for all $i\geq N$. 
\end{prop}

\begin{proof}
	We have that $(f_1-f_2)g = f_1g-f_2g = \sum_{i=0}^{N+m} (c_i-d_i)\t^i$. Since $c_i=d_i$ for all $i\geq N$, we have that $\deg_\t((f_1-f_2)g) < N$, so $\deg_\t(f_1-f_2)<N-\deg_\t(g) = N-m$. Since $f_1-f_2 = \sum_{i=0}^n (a_i-b_i)\t^i$, this gives $a_i=b_i$ for all $i\geq N-m$. 
\end{proof}

\begin{lemma} \th\label{DMcommutingsetprincipal} 
	Let $\a,B,C \in K\{\t\}$ be nonzero. Let $L$ be the perfect closure of $K$ and define the set $S_\a = \{x \in L\{\t\}\mid x\a = \a x\}$. Define $I = \{x \in K\{\t\}\mid B\a x = Cx\a \}$. If $S_\a \subset K\{\t\}$, then there exists a $k \in I$ such that $I = kS_\a$. 
\end{lemma}

\begin{proof}
	If $I = \{0\}$, then $I = 0\cdot S_\a$, so we may assume that $I \not=\{0\}$. This implies that $I$ contains a nonzero element $a$. Then $\deg_\t(B\a a) = \deg_\t(C a \a)$, which implies $\deg_\t(B) = \deg_\t(C)$. \\
	Let $k \in I$ be a nonzero element with the smallest $\t$-degree, which exists since we assumed that $I\not=\{0\}$. Let $s \in S_\a$, then $B\a ks = Ck\a s = Ck s \a$, as $k \in I$ and $s \in S_\a$. Since $S_\a \subset K\{\t\}$, we have $ks \in K\{\t\}$. This implies that $ks \in I$, so $kS_\a \subset I$. \\
	We will now prove the other inclusion. Let $a \in I$. Since $L$ is perfect, we have left and right division algorithms, see Theorem 3.1.13 and Exercise 3.1.8 of \cite{DM}, so there exist $b,r \in L\{\t\}$ such that $a = kb+r$ and $r=0$ or $\deg_\t(r)<\deg_\t(k)$. Since $L$ is perfect, there exists an element $\ol{b}$ such that $b\a = \a \ol{b}$.\\
	We have that $a \in I$, so $B\a a = Ca\a$, which implies
	$$B\a kb + B\a r = Ck b \a + Cr \a = Ck\a \ol{b} +Cr\a = B\a k \ol{b} + Cr \a,$$
	as $k \in I$. We have that $\deg_\t(B \a r) < \deg_\t(B \a k)$ and $\deg_\t(C \a r)< \deg_\t(C) + \deg_\t(\a k) = \deg_\t(B \a k)$. This means that $B\a kb$ and $B\a k \ol{b}$ have the same coefficient of $\t^i$ for all $i \geq \deg_\t(B \a k)$. By \th\ref{DMcoefpoly} $b$ and $\ol{b}$ have the same coefficient of $\t^i$ for all $i\geq \deg_\t(B\a k) - \deg_\t(B \a k) = 0$, so $b = \ol{b}$. \\
	This gives $b \in S_\a \subset K\{\t\}$, so $kb \in I$. Since $I$ is closed under addition, we get $r = \a-kb \in I$. We have $r=0$ or $\deg_\t(r)<\deg_\t(k)$ and $k$ was the nonzero element of smallest degree, so $r=0$. This gives $a = kb$ with $b \in S_\a$, so $I = kS_\a$ as needed.
\end{proof}

\begin{thm} \th\label{DMEnd=phiA} 
	Let $L$ be the perfect closure of $K$. Then $\End_K(\phi) = \End_L(\phi) = \phi(A) \subset K\{\t\}$. 
\end{thm}

\begin{proof}
	By Theorem 3.4.1 of \cite{DM} we get that $\End_K(\phi)$ and $\End_L(\phi)$ are free $A$-modules of rank at most $1$, where the multiplication is given by $a\cdot f = \phi_a f$. Since the identity is always a morphism, this gives $\End_K(\phi) = \End_L(\phi) = \phi(A)$. 
\end{proof}

\begin{lemma} \th\label{DMphiFeq} 
	We have $\phi(F) = \{x \in K(\t)\mid x\phi_T = \phi_T x\}$.
\end{lemma}

\begin{proof}
	Let $u \in \{x \in K(\t)\mid x\phi_T = \phi_Tx\}$ and write $u = \frac{a}{b}$ with $a,b \in K\{\t\}$. \\
	By the Ore condition there exist $B,C \in K\{\t\}\bs \{0\}$ such that $B(\phi_Tb) = C(b\phi_T)$. This gives 
	$$u\phi_T = \frac{a}{b}\cdot\frac{\phi_T}{1} = \frac{a\phi_T}{1\cdot b} = \frac{a\phi_T}{b} \text{ and } \phi_T u = \frac{\phi_T}{1}\cdot \frac{a}{b} = \frac{B\phi_T a}{Cb}.$$ Since $1\cdot(Cb) = C \cdot b$, the equality $u\phi_T = \phi_Tu$ gives $B\phi_Ta = Ca\phi_T$. \\
	Let $I = \{x \in K\{\t\}\mid B\phi_Tx = Cx\phi_T\}$ and let $L$ be the perfect closure of $K$. By \th\ref{DMEnd=phiA} we have that $\{s \in L\{\t\}\mid s\phi_T = \phi_Ts\} = \End_L(\phi) = \phi(A) \subset K\{\t\}$, so we may apply \th\ref{DMcommutingsetprincipal}. This gives that there exists a $k \in \phi(A)$ such that $I = k\phi(A)$. Since $a,b \in I$ per construction, there exist $a',b'\in \phi(A)$ such that $a = ka'$ and $b= kb'$. We find $\frac{a}{b} = \frac{ka'}{kb'} = \frac{a'}{b'} \in \phi(F)$, so $\{x \in K(\t) \mid x\phi_T = \phi_tx\} \subset \phi(F)$.\\
	For the converse, let $u \in \phi(F)$. Then we can write $u = \frac{a}{b}$ with $a,b \in \phi(A)$. Since $\End_K(\phi) = \phi(A)$, we get $a\phi_T = \phi_T a$ and $b\phi_T = \phi_T b$. This implies that 
	$$\phi_Tu = \frac{\phi_T}{1}\cdot \frac{a}{b} = \frac{\phi_Ta}{b\cdot 1} = \frac{a\phi_T}{b} = \frac{a}{b}\cdot \frac{\phi_T}{1} = u \phi_T.$$ 
	We conclude that $u \in \{x \in K(\t) \mid x\phi_T = \phi_Tx\}$, so $\phi(F) = \{x \in K(\t)\mid x\phi_T = \phi_Tx\}$.  
\end{proof}

\begin{assumption}
	We assume that the language on $K(\t)$ is $\L_{\t,\y}$. 
\end{assumption}

\begin{thm} \th\label{DMphiFdiophmap} 
	Give the language $\L_T$ to $F$. Then the map $\phi:F \ra K(\t)$ is an effective Diophantine map.
\end{thm}

\begin{proof}
	First note that the extended version of $\phi$ is still an injective map. We are going to use \th\ref{H10PeffectiveDioph}. 
	\begin{itemize} 
		\item By \th\ref{DMphiFeq} we know that 
		$$\phi(F) = \{u \in K(\t)\mid u \phi_T = \phi_Tu\} = \{u \in K(\t)\mid u(\y(T)+\t) = (\y(T)+\t)u\},$$
		which is a Diophantine set.
		\item We have that $\phi(0)= 0$, $\phi(1) = 1$ and $\phi(T) = \y(T) + \t$, 	which are all Diophantine expressions in the constants of $\L_{\t,\y}$.
		\item Let $* \in \{+,\cdot\}$. Since $\phi$ is a ring homomorphism we get
		\begin{align*}
			\phi(\{(a,b,a*b)\mid a,b \in F\}) &= \{(\phi(a),\phi(b),\phi(a)*\phi(b))\mid a,b \in F\} \\
			&= \{(x,y,z) \in K(\t)^3\mid x,y\in \phi(F),\ z = x*y\}.
		\end{align*}
		These are Diophantine sets since $\phi(F)$ is a Diophantine set.
	\end{itemize} 
	We have checked all the conditions of \ref{H10PeffectiveDioph}, so $\phi$ is an effective Diophantine map.
\end{proof}

\begin{lemma} \th\label{DMrecreprKtauDiv}
	If $L\{\t\}$ has a recursive representation $j: L\{\t\} \ra \N$, then $L(\t)$ has a recursive representation $J: L(\t) \ra \N^2$. 
\end{lemma} 

\begin{proof}
	The idea of the construction is the same as creating the recursive representation of $\Frac(R)$ from the recursive representation of a commutative integral domain $R$. We will do this proof informally.\\
	We have that $j$ is a recursive representation, so $j(L\{\t\})$ is a recursive set. This implies that $j(L\{\t\})^2$ is a recursive, so recursively enumerable set. Let $f$ be a function that injectively and recursively enumerates $j(L\{\t\})$ and set $(a_n,b_n) = f(n)$. \\
	We will first define a recursive function $H$ by strong recursion.\\
	We do a case-distinction to define $H(n+1)$.\\
	If $b_{n+1} = j(0)$, we set $H(n+1) = (j(0),j(0))$. This case is discarding division by zero. \\
	The next case corresponds to checking whether $\frac{j^{-1}(a_i)}{j^{-1}(b_i)} = \frac{j^{-1}(a_{n+1})}{j^{-1}(b_{n+1})}$. We can compute by Gauss-elimination the coefficients of some $A_i,B_i \in L\{\t\}$ such that $B_ij^{-1}(b_i) = A_ij^{-1}(b_{n+1})$. We can get the coefficients of elements of $L\{\t\}$ with a recursive function by \th\ref{DMrecreprKtau} and a Gauss-elimination uses only recursive case distinctions and field operations, so $B_i$ and $A_i$ can be determined with a recursive function. Then the check whether $\frac{j^{-1}(a_i)}{j^{-1}(b_i)} = \frac{j^{-1}(a_{n+1})}{j^{-1}(b_{n+1})}$ reduces to computing whether or not $B_ij^{-1}(a_i) = A_ij^{-1}(a_{n+1})$, which is a check in $L\{\t\}$, so it can be done with a recursive function. This means we can set $H(n+1) = H(i)$ if for some $0\leq i \leq n$ we have $\frac{j^{-1}(a_i)}{j^{-1}(b_i)} = \frac{j^{-1}(a_{n+1})}{j^{-1}(b_{n+1})}$. \\
	If $(a_{n+1},b_{n+1})$ does not correspond to a duplicate or to division by zero, we set $H(n+1) = (a_{n+1},b_{n+1})$. \\
	Now define $J(L(\t)) \ra \N^2$ by $J(r) = H((\mu n)(f(n) = (j(a),j(b))))$ if $r = \frac{a}{b}$ with $a,b \in L\{\t\}$. To check that $J$ is well-defined, let $r = \frac{a}{b} = \frac{a'}{b'}$ with $a,b,a',b' \in L\{\t\}$. Let $n,n'$ be minimal such that $f(n)= (j(a),j(b))$ and $f(n') =(j(a'),j(b'))$. Then by construction of $H$ we have that $H(n) = H(n')$, so $J$ is well-defined.\\
	Since $j$ is injective, we get by the construction of $H$ that $J$ is also injective.\\
	Next, we will show that $J(L(\t))$ is a recursive set. Let $(m,k) \in \N^2$. Since $j(L\{\t\})^2$ is a recursive set, we can check whether $(m,k) \in j(L\{\t\})^2$. Since $J(L(\t)) \subset H(\N) \subset j(L\{\t\})^2$, we have $(m,k)\notin J(L(\t))$ if $(m,k) \notin j(\L\{\t\})^2$. If $(m,k) \in j(L\{\t\})^2 = f(\N)$, there exists a unique $n$ such that $f(n) = (a_n,b_n) = (m,k)$. If $H(n) \not= (m,k)$, we get $(m,k) \notin J(L(\t))$, since $J(L(\t)) \subset H(\N)$. If $H(n) = (m,k)$ and $k\not= j(0)$, we get $(m,k) \in J(L(\t))$. If $H(n) = (m,k)$ and $k = j(0)$, we get that $(m,k)$ represents a division by zero fraction, so $(m,k) \notin J(L(\t))$. All these steps can be modelled by a recursive function, so $j(L(\t))$ is a recursive set.\\
	To get that the translates $T_+$ and $T_\ti$ of $+$ and $\cdot$ are recursive functions, we note that we can recursively get nonzero $b'$ and $d'$ such that $b'b = d'd$ and $a'$ and $d'$ such that $a'a=d'd$, respectively. This means that we can recursively get $x,y \in L\{\tau\}$ such that $\frac{x}{y}$ equals $\frac{a}{b}+\frac{c}{d}$ and $\frac{a}{b}\cdot \frac{c}{d}$, respectively. Then we can use the recursive function $H((\mu n)(f(n)=(j(x),j(y))))$ to get the correct representation of the fraction.\\
	Finally we have that $-\frac{a}{b} = \frac{-a}{b}$, as $\frac{a}{b} + \frac{-a}{b} = \frac{1\cdot a + 1\cdot (-a)}{1\cdot b} = \frac{0}{b} = 0$ and $1/\frac{a}{b} = \frac{b}{a}$, as $\frac{a}{b}\frac{b}{a} = \frac{1\cdot b}{1\cdot b}= 1$. Since we can recursively get a pair $(j(x),j(y))$ that represents $r = \frac{x}{y}$, also the translations of minus and division are recursive functions. 
\end{proof} 

\begin{thm} \th\label{DMKtauDivUndecidable} 
	Hilbert's tenth problem over $(K(\t),\L_{\t,\y})$ exists and has a negative answer.
\end{thm} 

\begin{proof}
	The ring $\F_p(\y(T))$ is a recursive ring, so by \th\ref{DMrecreprKtau,DMrecreprKtauDiv} we get that $\F_p(\y(T))(\t)$ has a recursive representation. Since $S_c(\L_{\t,\y}) \subset \F_p(\y(T))(\t) \subset S_c'(\L_{\t,\y})$, this gives that Hilbert's tenth problem over $(K(\t),\L_{\t,\y})$ exists.\\
	If $p\geq 3$, we get with \cite{PheidasFpT} that Hilbert's tenth problem over $(F,\L_T)$ has a negative answer. If $p=2$, we get the same result by \cite{H10PF2(T)}. By \th\ref{DMphiFdiophmap} we have that $\phi$ is an injective effective Diophantine map. This gives with \ref{H10PR2deci+deff=>R1deci} that Hilbert's tenth problem over $(K(\t),\L_{\t,\y})$ has a negative answer. 
\end{proof} 

\subsection{MRDP Holds in $\F_{q^n}\{\t\}$} \label{parDMMRDP} \phantom{=}\\
In \cite{FqtMRDP} Demeyer proves that MRDP holds in the ring $\F_q[T]$ if we add an extra constant to our language. We will transfer this result to the ring $\F_{q^n}\{\t\}$, where we still have that $a^q \t = \t a$ for all $a \in \F_{q^n}$.  

\begin{notation}
	For the rest of this section, let $\a$ be an element of the algebraic closure of $\F_q$ such that $\F_q(\a) = \F_{q^n}$ and let $\b$ be an element of the algebraic closure of $\F_p$ such that $\F_p(\b) = \F_q$. 
\end{notation} 

\begin{lemma} \th\label{DMphiArec}
	The function $\phi: A \ra \F_{q^n}(\y(T))\{\t\}$ is recursive.
\end{lemma} 

\begin{proof}
	By \th\ref{DMrecreprKtau}, $\F_{q^n}(\y(T))\{\t\}$ has a recursive representation, so it makes sense to speak of recursiveness of functions from $A$ to $\F_{q^n}(\y(T))\{\t\}$.\\
	Let $C: \N\ti A \ra A$ be the function such that $C(i,\sum_{j=0}^{n} a_jT^j) = a_i$ and let $\deg_T: A \ra \N$ be the degree function in $T$. Then $C$ and $\deg_T$ are recursive functions. \\ 
	By writing
	$$\y(x) = \sum_{i=0}^{\deg_T(x)} C(i,x) \y(T)^i,$$
	for $x\in A$, we conclude that $\y:A \ra \F_{q^n}(\y(T))\{\t\}$
	is a recursive function. Then
	$$\phi(x) = \sum_{i=0}^{\deg_T(x)} \y(C(i,x)) (\y(T) + \t)^i,$$
	is also a recursive function. 
\end{proof} 

For a ring $R$ with language $\L$ we can define a language $\L^{(k)}$ on $R^k$ that gives the same Diophantine sets as $(R,\L)$. This means that $S \subset R^{kl}$ is a Diophantine set seen via $(R,\L)$ if and only if it is a Diophantine set seen via $(R^k,\L^{(k)})$, see also Lemma 5.1.22 of \cite{ArtikelDioph}. 

\begin{defi} \th\label{DMdefcartlanguage}
	Let $(R,\L)$ be a set with an interpretation of the language $\L$. We define the language $\L^{(k)}$ and interpret it on $R^k$ as follows:
	\begin{itemize}
		\item For all $c \in \L$ a constant, $c \in \L^{(k)}$ is a constant and $c^{R^k} = (c^R,0^R,\ldots,0^R)$;
		\item For all $f \in \L$ an $r$-ary function, $f \in \L^{(k)}$ is an $r$-ary function and \\$f^{R^k}(\vec{x_1},\ldots,\vec{x_k}) = (f^R(x_{11},\ldots,x_{k1}),0^R,\ldots,0^R)$;
		\item For all $S \in \L$ an $r$-ary relation, $S \in \L^{(k)}$ is an $r$-ary relation and \\$S^{R^k} = \{((x_1,\vec{0}),\ldots,(x_r,\vec{0})) \in R^{kr} \mid \vec{x} \in S^R\}$;
		\item For $1\leq i \leq k$ the $1$-ary function $\pi_i$ is in $\L^{(k)}$ and $\pi_i(\vec{x}) = (x_i,0^R,\ldots,0^R)$.
	\end{itemize}
\end{defi}

\begin{thm} \th\label{DMMRDPFqntau}
	MRDP holds in the ring $\F_{q^n}\{\t\}$ with the language $\L_{\t,\y}\cup \{\a,\b\}$. 
\end{thm} 

\begin{proof}
	We have the relation $\t \a = \a^q \t$ and there exist a relation $\a^n = r(\a)$ with $n$ minimal and $r(X) \in \F_q[X]$ of degree at most $n-1$. This implies that we can write an element of $\F_{q^n}\{\t\}$ uniquely in the form $f_0+\a f_1+\ldots +\a^{n-1}f_{n-1}$ with $f_i \in \F_q\{\t\}$. Using this and the Carlitz module $\phi$, we can define a map $d: (\F_q[T])^n \ra \F_{q^n}\{\t\}$ by 
	$$d(f_0,\ldots,f_{n-1}) = \phi(f_0)+\a \phi(f_1)+\ldots \a^{n-1}\phi(f_{n-1}).$$
	We have that $\deg_T(f)$ and $\deg_\t(\phi(f))$ are equal for all $f \in \F_q[T]$. Furthermore $\F_q\{\t\}$ and $\F_q[T]$ both have $q^{d+1}$ elements of degree $d$. Since $\phi$ is injective, this gives that $\phi$ is a bijection between $\F_q[T]$ and $\F_q\{\t\}$. Using the fact that we can write each element of $\F_{q^n}\{\t\}$ uniquely in the form $f_0+\a f_1+\ldots + \a^{n-1} f_{n-1}$, we get that $d$ is bijective. \\
	We have that $d$ is given by a Diophantine expression in $\phi$, for the exact definition see Definition 5.3.1 of \cite{ArtikelDioph}, namely
	$$d(f_0,\ldots,f_{n-1}) = y \LRa y = \phi(f_0)+\a \phi(f_1)+\ldots + \a^{n-1}\phi(f_{n-1}).$$
	With \th\ref{DMphiAdioph} we have that $\phi$ is a Diophantine map if we give $A$ the language $\L_T$. Since $\phi(\b) = \b$ is a constant of $\L_{t,\y}\cup\{\a,\b\}$ we get that $\phi: (A,\L_{T}\cup\{\b\})\ra (K\{\t\},\L_{\t,\y}\cup\{\a,\b\})$ is also a Diophantine map. With Proposition 5.3.4 of \cite{ArtikelDioph} this gives that $d$ is a Diophantine map. \\
	We have that $d$ is the composition of recursive functions, since $\phi$ is recursive by \th\ref{DMphiArec}, so $d$ is a recursive function.\\
	In \cite{FqtMRDP} Demeyer proves that MRDP holds in $(\F_q[T],\L_T\cup \{\b\})$. With Proposition 6.2.13 of \cite{ArtikelDioph} this gives that MRDP holds in $(\F_q[T]^n,(\L_T\cup\{\b\})^{(n)})$. Then with Theorem 7.5.1 of \cite{ArtikelDioph} we get that MRDP holds in $(\F_{q^n}\{\t\},\L_{\t,\y}\cup \{\a,\b\})$, since $d$ is a bijective Diophantine and recursive map. 
\end{proof} 

\subsection{Reductions for $\F_{q^n}\pownon{\t}$} \label{parDMpower} \phantom{=}\\
In this paragraph, we prove that Hilbert's tenth problem over $\F_{q^n}\pownon{\t}$ 
has a positive respectively negative answer if Hilbert's tenth problem over the commutative ring $\F_p\pow{T}$ has a positive respectively negative answer. In the next paragraph, we do the same for the rings $\F_{q^n}\lau{\t}$ and $\F_p\lau{T}$. \\
In \cite{Fq[[T]]ADF}, Anscombe, Dittmann and Fehm prove that Hilbert's tenth problem over $(\F_q\lau{T},\L_T\cup\{\ord_T\})$ has a positive answer under the following condition:
\begin{conjecture}[Weak resolution of singularities] \th\label{DMWeakResOfSing} 
	\nieuw{For all fields $L_0$ and all nontrivial finitely generated extensions $L_1$ of $L_0$ such that there exists a valuation $v$ on $L_1$ with residue field $L_0$, there also exists a nontrivial discrete valuation $w:L_1 \ra \Z$ with residue field $L_0$.}
\end{conjecture} 
This result then implies that Hilbert's tenth problem over $(\F_q\lau{T},\L_T)$ and $(\F_q\pow{T},\L_T)$ have a positive answer under the same hypothesis. Since we do not need $\ord_T$, we discard it. This is allowed, since we have a positive answer and discarding parts of the language makes Hilbert's tenth problem easier. At this moment there are no unconditional results for a positive or negative answer to Hilbert's tenth problem over these rings with the language $\L_T$. 

For the proof, we need the following lemma. This is a variant from a lemma of \cite{Ct1t2PZSurvey}.
\begin{lemma} \th\label{Ct1t2Addalpha} 
	Let $(R,\L)$ be an integral domain with a language $\L$, such that $\L$ is the language of rings with possibly more constants and relations. Suppose that Hilbert's tenth problem over $(R,\L)$ has a positive answer. Let $K$ be the fraction field of $R$ and let $\ol{K}$ be the algebraic closure of $K$. Let $\a \in \ol{K}$ be integral over $R$ such that its minimal equation has coefficients in $R_0 = S_c'(\L)$. Define the language $\L_\a = \L \cup \{c_{\a},\inR\}$ where $c_\a$ is a constant and $\inR$ a unary relation. We interpret $\L_\a$ in $R[\a]$ as follows:
	\begin{itemize}
		\item $+,\cdot$ are the usual addition and multiplication;
		\item $c^{R[\a]} = c^R$ for all constants $c \in \L$;
		\item $S^{R[\a]} = S^R$ for all relations $S \in \L$;
		\item $c_\a^{R[\a]} = \a$;
		\item $\inR(x)$ if and only if $x \in R$.
	\end{itemize}
	Then Hilbert's tenth problem over $(R[\a],\L_\a)$ has a positive answer.
\end{lemma}

\begin{proof}
	An arbitrary Diophantine statement in $\L_\a$ is of the form 
	\begin{align} \label{Ct1t2EqStart} 
		\ex \vec{x} \in R[\a]^l, \vec{y} \in R[\a]^k,\ \va 1\leq i \leq m,\ f_i(\vec{x},\vec{y})=0,\ \va 1\leq i \leq l,\ \inR(x_i), \notag \\ 
		\va 1\leq i \leq m',\ (g_{i1}(\vec{x},\vec{y}),\ldots,g_{in_i}(\vec{x},\vec{y}) \in S_i,
	\end{align} 
	where $f_i,g_{ij}$ are polynomials with coefficients in $S_c'(\L_\a)$ and $S_i$ are relations of $\L$. \\
	Since $\a$ is integral over $R$, there exists a smallest $n$ such that $1,\a,\ldots,\a^n$ generate $R[\a]$ as a module over $R$. Define new variables $y_{i0},\ldots,y_{in}$ for $1\leq i \leq k$ and set $y_i = y_{i0} + \ldots + y_{in}\a^n$. We can substitute this into the equations $f_i(\vec{x},\vec{y})$. Then we can use that \nieuw{$R_0[\a]$ is generated by $\{1,\ldots,\a^{n}\}$ over $R_0$} to rewrite these equations into the form
	\begin{align}
		\label{Ct1t2Eq1} f_{i0}(\vec{x},\vec{y_1},\ldots,\vec{y_k}) + \ldots + \a^n f_{in}(\vec{x},\vec{y_1},\ldots,\vec{y_k}) = 0,\\
		(g_{i10}(\vec{x},\vec{y_1},\ldots,\vec{y_k}) + \ldots + \a^n  g_{i1n}(\vec{x},\vec{y_1},\ldots,\vec{y_k}),\ldots ,\phantom{=========} \notag \\ g_{in_i0}(\vec{x},\vec{y_1},\ldots,\vec{y_k}) + \ldots + \a^n g_{in_in}(\vec{x},\vec{y_1},\ldots,\vec{y_k})) \in S_i, \label{Ct1t2Eq2}
	\end{align} 
	where the $f_{ij}$ and $g_{iji'}$ are polynomials with coefficients in \nieuw{$R_0 = S_c'(\L)$}. \\
	Since the set $\{1,\ldots,\a_n\}$ is linearly independent over $K \supset R$, we get that Equation \ref{Ct1t2Eq1} is equivalent to
	\begin{align} \label{Ct1t2Eq3}  
		\va 0\leq j \leq n,\ f_{ij}(\vec{x},\vec{y_1},\ldots,\vec{y_k}) = 0.
	\end{align} 
	To rewrite Equation \ref{Ct1t2Eq2} we note that $S_i \subset R^{n_i}$ and that $$g_{ij0}(\vec{x},\vec{y_1},\ldots,\vec{y_k}) + \ldots + \a^n g_{ijn}(\vec{x},\vec{y_1},\ldots,\vec{y_k}) \in R$$
	is equivalent to 
	$$\va 1\leq i' \leq n,\ g_{iji'}(\vec{x},\vec{y_1},\ldots,\vec{y_k}) = 0.$$ 
	This gives that Equation \ref{Ct1t2Eq2} is equivalent to 
	\begin{align} \label{Ct1t2Eq4} 
		(g_{i10}(\vec{x},\vec{y_1},\ldots,\vec{y_k}),\ldots,g_{in_i0}(\vec{x},\vec{y_1},\ldots,\vec{y_k})) \in S_i\\
		\va 1\leq i' \leq n,\ g_{iji'}(\vec{x},\vec{y_1},\ldots,\vec{y_k}) = 0.
	\end{align} 
	Combining Equations \ref{Ct1t2Eq3} and \ref{Ct1t2Eq4} gives that our arbitrary Diophantine Equation \ref{Ct1t2EqStart} is equivalent to 
	\begin{align*} 
		\ex \vec{x} \in R^l, \ex \vec{y_1},\ldots,\vec{y_k} \in R^n,\\ 
		\va 1\leq i \leq m,\ \va 0\leq j \leq n,\ f_{ij}(\vec{x},\vec{y_1},\ldots,\vec{y_k}) = 0\\
		\va 1\leq i \leq m',(g_{i10}(\vec{x},\vec{y_1},\ldots,\vec{y_k}),\ldots,g_{in_i0}(\vec{x},\vec{y_1},\ldots,\vec{y_k})) \in S_i\\
		\va 1\leq i' \leq n,\ g_{iji'}(\vec{x},\vec{y_1},\ldots,\vec{y_k}) = 0.
	\end{align*} 
	This is a Diophantine statement in $(R,\L)$, so by the positive answer to Hilbert's tenth problem over $(R,\L)$, we get that Hilbert's tenth problem over $(R[\a],\L_\a)$ has a positive answer.
\end{proof}

\begin{lemma} \th\label{DMaddconstFq[[T]]}
	Suppose that Hilbert's tenth problem over $(\F_p\pow{T},\L_T)$ has a positive answer. Then Hilbert's tenth problem over $(\F_q\pow{T},\L_T\cup\{\b\})$ has a positive answer.
\end{lemma} 

\begin{proof}
	This follows directly by applying \th\ref{Ct1t2Addalpha} on the ring $\F_p\pow{T}$ and the constant $\b$. 
\end{proof} 

\begin{lemma} \th\label{DMaddrelFq[[T]]}
	Suppose that Hilbert's tenth problem over $(\F_q\pow{T},\L_T\cup\{\b\})$ has a positive answer. Then for a fixed $m \in \Z_{\geq 1}$, Hilbert's tenth problem over $(\F_q\pow{T},\L_T\cup \{\b,\TextIn\F_q\pow{T^m}\})$ has a positive answer.
\end{lemma} 

\begin{proof}
	Let $S$ be such that $S^m = T$. Then $S$ is integral over $\F_q\pow{T}$. We are going to use \th\ref{Ct1t2Addalpha} with $\a = S$ and $R = \F_q\pow{T}$. This gives that Hilbert's tenth problem over $(\F_q\pow{T}[S],\L_T\cup\{\b,S,\TextIn\F_q\pow{T}\})$ has a positive answer. We have that 
	\begin{align*}
		\F_q\pow{T}[S] &= \{\sum_{i=0}^\infty a_{0i}T^i + S\sum_{i=0}^\infty a_{1i}T^i + \ldots + S^{m-1}\sum_{i=0}^\infty a_{m-1,i}T^i \mid a_{ji} \in \F_q\}\\ 
		&= \{\sum_{i=0}^\infty \sum_{k=0}^{m-1} a_{ki}S^{k+mi}\mid a_{ji} \in \F_q\} = \F_q\pow{S}.
	\end{align*}
	This gives that Hilbert's tenth problem over $(\F_q\pow{S},\L_R\cup\{\b,S^m,S,\TextIn\F_q\pow{S^m}\})$ has a positive answer. By replacing the symbol $S$ by the symbol $T$ and by noticing that the constant $S^m$ is unnecessary, we get that Hilbert's tenth problem over $(\F_q\pow{T},\L_T\cup \{\b,\TextIn\F_q\pow{T^m}\}$ has a positive answer. 
\end{proof} 

\begin{defi}
	Since we have that $\t a = a^q \t$ for all $a \in \F_{q^n}$ and $\a^n = r(\a)$ for some $r(X) \in \F_q[X]$ we can write each element of $\F_q[\a]\pownon{\t}$ uniquely in the form $f_0+\a f_1+\ldots + \a^{n-1} f_{n-1}$ with $f_i \in \F_q\pownon{\t}$. Define the map $d: \F_q[\a]\pownon{\t} \ra \F_q\pow{T}^n$ by
	$$d(f_0(\t)+\a f_1(\t)+\ldots \a^{n-1}f_{n-1}(\t)) = (f_0(T),\ldots,f_{n-1}(T)).$$ 
\end{defi} 

\begin{assumption}
	We will assume that the language on $\F_q[\a]\pownon{\t}$ is $\L_R\cup \{\t,\a,\b\}$ and that the language on $\F_q\pow{T}^n$ is $(\L_T\cup\{\b,\TextIn\F_q[[T^n]]\})^{(n)}$. 
\end{assumption}

\begin{lemma} \th\label{DMpowerseriesconstsetsadd}
	The following sets are Diophantine:
	\begin{enumerate}
		\item $\{d(0)\},\ \{d(1)\},\ \{d(\t)\}$, $\{d(\a)\}$ and $\{d(\b)\}$;
		\item $d(\F_q[\a]\pownon{\t})$;
		\item $d(\{(f,g,f+g)\mid f,g \in \F_q[\a]\pownon{\t}\})$.
	\end{enumerate}
\end{lemma}

\begin{proof} \phantom{=}
	\begin{enumerate}
		\item We have that $\{d(0)\} = \{(0,\ldots,0)\}$, $\{d(1)\} = \{(1,0,\ldots,0)\}$, $\{d(\t)\} = \{(T,0,\ldots,0)\}$ and $\{d(\b)\} = \{(\b,0,\ldots,0)\}$ are all Diophantine sets, since the elements $(0,\ldots,0)$, $(1,0,\ldots,0)$, $(T,0,\ldots,0)$ and $(\b,0,\ldots,0)$ are constants of $\L_T^{(n)}$. Finally, we have that $$\{d(\a)\} = \{(0,1,0,\ldots,0)\} = \{x \in \F_q\pow{T}^n\mid \pi_1(x) = \vec{0},\ \pi_2(x) = (1,0,\ldots,0),\ \va 3\leq i \leq n,\ \pi_i(x)=\vec{0}\}$$ 
		is also a Diophantine set.
		\item We have that $d(\F_q[\a]\pownon{\t}) = \F_q\pow{T}^n$, which is a Diophantine set.
		\item We have 
		\begin{align*} 
			&d(\{f,g,f+g,\mid f,g \in \F_q[\a]\pownon{\t}\})  =\\
			&d(\{(f_0+\ldots + \a^{n-1} f_{n-1},g_0+\ldots + \a^{n-1}g_{n-1},(f_0+g_0)+\ldots+\a^{n-1}(f_{n-1}+g_{n-1}))\mid f_i,g_i \in \F_q\pownon{\t}\})\\
			&= \{((f_0,\ldots,f_{n-1}),(g_0,\ldots,g_{n-1}),(f_0+g_0,\ldots,f_{n-1}+g_{n-1}))\mid f_i,g_i \in \F_q\pow{T}\}\\
			&= \{(\vec{f},\vec{g},\vec{h})\mid \vec{f},\vec{g} \in \F_q\pow{T}^n,\ \vec{h} = \vec{f}+\vec{g}\}.
		\end{align*} 
		This gives that $d(\{(f,g,f+g)\mid f,g \in \F_q[\a]\pownon{\t}\})$ is a Diophantine set.
	\end{enumerate}
\end{proof}

\begin{lemma} \th\label{DMpowerseriesmult} 
	The set $d(\{(f,g,f\cdot g)\mid f,g \in \F_q[\a]\pownon{\t}\})$ is Diophantine.
\end{lemma} 

\begin{proof}
	We have that $\a$ is a nonzero element of $\F_{q^n}$, so $\a^{q^n-1} = 1$. Since $q^n\equiv 1 \mod q^n-1$ we get $q^{ni} \equiv 1 \mod q^n-1$ for all $i \in \N$. This gives 
	$$\a^{q^{ni}}=\a^{(q^n-1)c+1} = (\a^{q^n-1})^c \cdot \a = 1^c\cdot \a = \a,$$ so that $\a^{q^{ni+j}} = (\a^{q^{ni}})^{q^j} = \a^{q^j}$.\\
	Now let $f_0+\ldots + \a^{n-1}f_{n-1}$ and $g_0+\ldots + \a^{n-1}g_{n-1}$ be elements of $\F_q[\a]\pownon{\t}$ in standard form. Write $f_k = \sum_{i=0}^\infty a_{ki}\t^i$ with $a_{ki} \in \F_q$ and define $f_{kj} = \sum_{i=0}^\infty a_{k,ni+j}\t^{ni+j}$. Then
	\begin{align*}
		&\phantom{=} (f_0+\ldots+\a^{n-1}f_{n-1})(g_0+\ldots + \a^{n-1}g_{n-1}) \\
		&= \sum_{k=0}^{n-1} \sum_{l=0}^{n-1} \a^k \left(\sum_{i=0}^\infty a_{ki}\t^i \right) \a^l g_l \\
		&= \sum_{k=0}^{n-1} \sum_{l=0}^{n-1} \a^k \left(\sum_{i=0}^\infty a_{ki}(\a^{q^i})^l \t^i\right) g_l \\
		&= \sum_{k=0}^{n-1} \sum_{l=0}^{n-1} \a^k \left(\sum_{j=0}^{n-1} \sum_{i=0}^\infty \left(\a^{q^{ni+j}}\right)^l a_{k,ni+j}\t^{ni+j}\right) g_l \\
		&= \sum_{k=0}^{n-1} \sum_{l=0}^{n-1} \left(\sum_{j=0}^{n-1} \a^k \sum_{i=0}^\infty \left(\a^{q^{j}}\right)^l a_{k,ni+j}\t^{ni+j}\right) g_l \\
		&= \sum_{k=0}^{n-1} \sum_{l=0}^{n-1} \sum_{j=0}^{n-1} \a^k (\a^{q^j})^l f_{kj}g_l.
	\end{align*}  
	Since $n$ is fixed and we have the relation $\a^n = r(\a)$, there exist polynomials $P_i$ \nieuw{with coefficients in $\F_q$} such that
	\begin{align*} 
		&\sum_{k=0}^{n-1} \sum_{l=0}^{n-1} \sum_{j=0}^{n-1} \a^k (\a^{q^j})^l f_{kj}g_l \\
		&= P_0(f_{00},\ldots,f_{n-1,n-1},g_0,\ldots,g_{n-1}) + \ldots + \a^{n-1}P_{n-1}(f_{00},\ldots,f_{n-1,n-1},g_0,\ldots,g_{n-1}).
	\end{align*} 
	This implies that
	\begin{align*}
		&\phantom{=}d(\{(f,g,f\cdot g)\mid f,g \in \F_q[\a]\pownon{\t}\}) \\
		&= \{(\vec{f},\vec{g},\vec{h}) \mid \ex f_{kj}'\in \F_q[[T^n]], \ex f_{kj} \in \F_q\pow{T},\ \sum_{j=0}^{n-1} f_{kj} = f_k,\ f_{kj} = T^j\cdot f_{kj}',\\
		&\phantom{=========} h_i = P_i(f_{00},\ldots,f_{n-1,n-1},g_0,\ldots,g_{n-1}) \}.
	\end{align*} 
	Since we have the relation $\TextIn\F_q[[T^n]]$ and the constant $\b$, we get that $d(\{(f,g,f\cdot g)\mid f,g, \in \F_q[\a]\pownon{\t}\})$ is a Diophantine set.
\end{proof} 

\begin{thm} \th\label{DMpowerseriesdecidable}
	The map $d$ is an effective Diophantine map. If Hilbert's tenth problem over $(\F_p\pow{T},\L_T)$ has a positive answer, then Hilbert's tenth problem over $(\F_{q^n}\pownon{\t},\L_{\t,\y}\cup\{\a,\b\})$ has a positive answer. 
\end{thm} 

\begin{proof}
	We use \th\ref{H10PeffectiveDioph}. With \th\ref{DMpowerseriesconstsetsadd,DMpowerseriesmult} we get the explicit Diophantine formulas for the required sets, so $d$ is an effective Diophantine map.\\ 
	We have that $\F_p[T]$ and $\F_q[T]$ and $\F_{q^n}$ are recursive rings and by \th\ref{DMrecreprKtau} it follows that $\F_{q^n}\{\t\}$ has a recursive representation. We have $S_c(\L_T) \subset \F_p[T]\subset S_c'(\L_T)$, $S_c(\L_T\cup\{\b\}) \subset \F_q[T] \subset S_c'(\L_T\cup\{\b\})$ and $S_c(\L_{\t,\y}\cup \{\a\}) \subset \F_{q^n}\{\t\} \subset S_c'(\L_{\t,\y}\cup \{\a\})$, which gives that Hilbert's tenth problem exists over the rings involved. \\ 
	By \th\ref{DMaddconstFq[[T]]} we get a positive answer over $(\F_q\pow{T},\L_T\cup\{\b\})$, if Hilbert's tenth problem over $(\F_p\pow{T},\L_T)$ has a positive answer. Using \th\ref{DMaddrelFq[[T]]} we can extend the language to $\L_T\cup \{\TextIn\F_q[[T^n]]\}$. It follows that Hilbert's tenth problem over $(\F_q\pow{T}^n,(\L_T\cup\{\TextIn\F_q[[T^n]]\})^{(n)})$ has a positive answer. Finally from \th\ref{H10PR2deci+deff=>R1deci} we get that Hilbert's tenth problem over $\F_{q^n}\pownon{\t}$ with the language $\L_{\t,\y}\cup\{\a\}$ has a positive answer. 
\end{proof} 

\begin{cor}
	If \th\ref{DMWeakResOfSing}, the weak resolution of singularities in characteristic $p$, holds, then Hilbert's tenth problem over $\F_{q^n}\pownon{\t}$ has a positive answer.
\end{cor} 

\begin{proof}
	\nieuw{In \cite{Fq[[T]]ADF} it is proven that if \th\ref{DMWeakResOfSing} holds, then Hilbert's tenth problem over $(\F_p\pow{T},\L_T)$ has a positive answer. With \th\ref{DMpowerseriesdecidable} this gives that under this hypothesis Hilbert's tenth problem over $\F_{q^n}\pownon{\t}$ has a positive answer.} 
\end{proof} 

\begin{thm} The Carlitz module can be extended to a map $\phi: \F_q\pow{T} \ra K\pownon{\t}$. If we give $\F_q\pow{T}$ the language $\L_T$ and $K$ the special $A$-field structure $\y_0$ from \th\ref{DMdefy0}, then this is an effective Diophantine map. \\
	\nieuw{If Hilbert's tenth problem over $\F_q\pow{T}$ has a negative answer, then Hilbert's tenth problem over $K\pownon{\t}$ has a negative answer.}
\end{thm} 

\begin{proof}
	Since $\phi$ is injective, we are going to check all conditions of \th\ref{H10PeffectiveDioph}. \\
	We have that 
	$$\t \left(\sum_{i=0}^\infty a_i\t^i\right) = \sum_{i=0}^\infty a_i^q \t^{i+1} = \left(\sum_{i=0}^\infty a_i^q \t^i\right) \t.$$ 
	Since $a^q = a$ only holds for elements of $\F_q \subset K$ we get 
	$$\phi(\F_q\pow{T}) = \F_q\pownon{\t} = \{x \in K\pownon{\t}\mid x\t = \t x\}.$$
	This means that $\phi(\F_q\pow{T})$ is a Diophantine set. \\
	We have that $\phi(T) = \t$, $\phi(1) = 1$ and $\phi(0) = 0$, which are all constants in the language of $K\pownon{\t}$. \\
	Since $\phi$ is still a ring homomorphism and $\phi(\F_q\pow{T})$ is a Diophantine set, we get that the sets 
	$$\{(\phi(a),\phi(b),\phi(a*b))\mid a,b\in \F_q\pow{T}\}$$ with $* \in \{+,\cdot\}$ are Diophantine sets.\\
	These are all conditions of \th\ref{H10PeffectiveDioph}, so $d$ is a Diophantine map.\\
	The rings $\F_p[T]$ and $\F_p$ are recursive rings and by \th\ref{DMrecreprKtau} also $\F_p\{\t\}$ has a recursive representation. We have the inclusions $S_c(\L_T) \subset \F_p[T] \subset S_c'(\L_T)$ and $S_c(\L_{\t,\y_0})\subset \F_p\{\t\} \subset S_c'(\L_{\t,\y_0})$, which gives that Hilbert's tenth problem exists over $(\F_q\pow{T},\L_T)$ and $(K\pownon{\t},\L_{\t,\y_0})$. Using \th\ref{H10PR2deci+deff=>R1deci} and contraposition we get that if Hilbert's tenth problem over $\F_q\pow{T}$ has a negative answer, then Hilbert's tenth problem over $K\pownon{\t}$ has a negative answer.
\end{proof} 

\subsection{\nieuw{Reductions for $\F_{q^n}\lau{\t}$}} \label{parDMlaurent} 
\begin{lemma} \th\label{DMaddconstFq((T))}
	\nieuw{If Hilbert's tenth problem over $(\F_p\lau{T},\L_T)$ has a positive answer, then Hilbert's tenth problem over $(\F_q\lau{T},\L_T\cup\{\b\})$ has a positive answer.}
\end{lemma}

\begin{proof}
	\nieuw{This follows immediately from \th\ref{Ct1t2Addalpha} with the ring $\F_p\lau{T}$ and constant $\b$.} 
\end{proof} 

\begin{lemma} \th\label{DMlaurentseriesAddRel}
	Fix some $m \in \Z_{>1}$. If Hilbert's tenth problem over \nieuw{$(\F_q\lau{T},\L_T\cup\{\b\})$} has a positive answer, then Hilbert's tenth problem over \nieuw{$(\F_q\lau{T},\L_T\cup\{\b,\TextIn\F_q\lau{T^m}\})$} has a positive answer.
\end{lemma} 

\begin{proof}
	We have that $\sqrt[m]{T}$ is integral over the commutative ring $\F_q\lau{T}$ and $\L_T\cup\{\b\}$ contains no other functions than addition and multiplication. This means we can use \th\ref{Ct1t2Addalpha} with $R = \F_q\lau{T}$ and $\a = \sqrt[m]{T}$. This gives that if Hilbert's tenth problem over $(\F_q\lau{T},\L_T\cup\{\b\})$ has a positive answer, then Hilbert's tenth problem over $(\F_q\lau{T}[\sqrt[m]{T}],\L_T\cup\{\b,\TextIn\F_q\lau{T},\sqrt[m]{T}\})$ has a positive answer.\\
	If we substitute $S = \sqrt[m]{T}$ we get $\F_q\lau{T}[\sqrt[m]{T}] = \F_q\lau{S^m}[S] = \F_q\lau{S}$. The language then becomes $\L_R \cup \{\b,\TextIn\F_q\lau{S^m},S,S^m\}$. If we remove $S^m$ from the language and do a further change of variables $S= T$, we get the result. 
\end{proof} 

\begin{assumption}
	We assume that the language on $\F_q\lau{T}^n$ is $(\L_T \cup \{\b,\TextIn\F_q((T^n))\})^{(n)}$ and the language on $\F_{q^n}\lau{\t}$ is \nieuw{$\L_{\t,\y_0}\cup\{\a,\b\}$}. 
\end{assumption} 

\begin{lemma}
	The map $d: \F_{q^n}\pownon{\t} \ra \F_q\pow{T}^n$ extends to a map $d':\F_{q^n}\lau{\t} \ra \F_q\lau{T}^n$. 
\end{lemma} 

\begin{proof}
	We have that $a^q \t = \t a$ for all $a \in \F_{q^n}$. If $a\not=0$ we can take inverses, which gives the equality $\t^{-1} a^{-q} = a^{-1} \t^{-1}$. We have that $\F_{q^n}$ is a finite field, so it is perfect. Since $q$ is a power of $p$, we can take $q$-th roots. Now take $a = \sqrt[q]{b^{-1}}$, then $b^{q^{-1}}t^{-1} = \t^{-1} b$ for all nonzero $b \in \F_{q^n}$. Using this and the relation $\a^n = r(\a)$ with $r(X) \in \F_q[X]$, we can write an element of $\F_q[\a]\lau{\t}$ uniquely in the form $f_0+\a f_1 +\ldots + \a^{n-1}f_{n-1}$ with $f_i \in \F_q\lau{\t}$. With this we define
	\begin{equation*} 
		d'(f_0+\a f_1+\ldots + \a^{n-1}f_{n-1}) = (f_0,\ldots,f_{n-1}). \qedhere \end{equation*} 
\end{proof} 

\begin{lemma} \th\label{DMlaurentseriesmult} 
	The set $d'(\{(f,g,f\cdot g)\mid f,g \in \F_{q^n}\lau{\t}\})$ is Diophantine.
\end{lemma} 

\begin{proof}
	We will argue in the same way as in the proof of \th\ref{DMpowerseriesmult}. \\
	It still holds that $\a^{q^{ni+j}} = \a^{q^j}$ for all $i \in \Z_{\geq 0}$. Since $\F_{q^n}$ is perfect and $q$ is a $p$-th power, we get from the equality $\a^{q^{ni}} = \a$ by taking a $q^{ni}$-th root that $\a = \a^{q^{-ni}}$. This gives that $\a^{q^{ni+j}} = \a^{q^j}$ for all $i \in \Z$.\\
	Take $f_0+\ldots + \a^{n-1}f_{n-1}$ and $g_0+\ldots + \a^{n-1}g_{n-1}$ elements of $\F_q[\a]\lau{\t}$ in the standard form. We write $f_k = \sum_{i=-N_{f_k}}^\infty a_{ki}\t^i$ with $a_{ki} \in \F_q$ and define $f_{kj} = \sum_{i=-N_{f_{kj}}}^\infty a_{k,ni+j}\t^{ni+j}$. Then
	\begin{align*}
		&\phantom{=} (f_0+\ldots+\a^{n-1}f_{n-1})(g_0+\ldots + \a^{n-1}g_{n-1}) \\
		&= \sum_{k=0}^{n-1} \sum_{l=0}^{n-1} \a^k \left(\sum_{i=-N_{f_k}}^\infty a_{ki}\t^i \right) \a^l g_l \\
		&= \sum_{k=0}^{n-1} \sum_{l=0}^{n-1} \a^k \left(\sum_{i=-N_{f_k}}^\infty a_{ki}(\a^{q^i})^l \t^i\right) g_l \\
		&= \sum_{k=0}^{n-1} \sum_{l=0}^{n-1} \a^k \left(\sum_{j=0}^{n-1} \sum_{i=-N_{f_{kj}}}^\infty \left(\a^{q^{ni+j}}\right)^l a_{k,ni+j}\t^{ni+j}\right) g_l \\
		&= \sum_{k=0}^{n-1} \sum_{l=0}^{n-1} \left(\sum_{j=0}^{n-1} \a^k \sum_{i=-N_{f_{kj}}}^\infty \left(\a^{q^{j}}\right)^l a_{k,ni+j}\t^{ni+j}\right) g_l \\
		&= \sum_{k=0}^{n-1} \sum_{l=0}^{n-1} \sum_{j=0}^{n-1} \a^k (\a^{q^j})^l f_{kj}g_l.
	\end{align*}  
	Since $n$ is fixed and we have the relation $\a^n = r(\a)$, there exist polynomials $Q_i$ \nieuw{with coefficients in $\F_q$} such that
	\begin{align*} 
		&\sum_{k=0}^{n-1} \sum_{l=0}^{n-1} \sum_{j=0}^{n-1} \a^k (\a^{q^j})^l f_{kj}g_l \\
		&= Q_0(f_{00},\ldots,f_{n-1,n-1},g_0,\ldots,g_{n-1}) + \ldots + \a^{n-1}Q_{n-1}(f_{00},\ldots,f_{n-1,n-1},g_0,\ldots,g_{n-1}).
	\end{align*} 
	This implies that
	\begin{align*}
		&\phantom{=}d(\{(f,g,f\cdot g)\mid f,g \in \F_q[\a]\lau{\t} \}) \\
		&\phantom{=}= \{(\vec{f},\vec{g},\vec{h}) \mid \ex f_{kj}'\in \F_q\lau{T^n}, \ex f_{kj} \in \F_q\lau{T},\ \sum_{j=0}^{n-1} f_{kj} = f_k,\\ 
		&\phantom{=========}f_{kj} = T^j\cdot f_{kj}',\ h_i = Q_i(f_{00},\ldots,f_{n-1,n-1},g_0,\ldots,g_{n-1}) \}.
	\end{align*} 
	By the construction of our language we get that $d(\{(f,g,f\cdot g)\mid f,g \in \F_q[\a]\pownon{\t}\})$ is a Diophantine set.
\end{proof} 

\begin{thm} \th\label{DMlaurentseriesdecidable}
	The map $d': \F_{q^n}\lau{\t}\ra \F_q\lau{T}^n$ is an effective Diophantine map.\\
	If Hilbert's tenth problem over $(\F_p\lau{T},\L_T)$ has a positive answer, then Hilbert's tenth problem over $\F_{q^n}\lau{\t}$ has a positive answer.
\end{thm} 

\begin{proof}
	We have that $d'$ is injective, so we are going to use \th\ref{H10PeffectiveDioph}.\\ 
	We have that $\{d'(0)\} = \{\vec{x} \in \F_q\lau{T}^n\mid \vec{x} = \vec{0}\}$, $\{d'(1)\} = \{(1,0,\ldots,0)\}$, $\{d'(\y_0(T))\} = \{d'(0)\} = \{\vec{0}\}$, \nieuw{$\{d'(\b)\} = \{(\b,0,\ldots,0)\}$} and
	$$\{d'(\a)\} = \{(0,1,0,\ldots,0)\} = \{\vec{x} \in \F_q\lau{T}^n\mid \pi_1(\vec{x})=0,\ \pi_2(\vec{x})=1,\va 3\leq i \leq n,\ \pi_i(\vec{x})=0\},$$ which are all Diophantine sets in the required form.\\
	We have $d'(\F_{q^n}\lau{\t}) = \F_q\lau{T}^n$, which is clearly a Diophantine set.\\
	We have that 
	\begin{align*}
		d'(\{(a,b,a+b)\mid a,b \in \F_{q^n}\lau{\t}\}) = \{(\vec{a},\vec{b},\vec{c})\mid \vec{a},\vec{b} \in \F_q\lau{T}^n,\ \vec{c} = \vec{a}+\vec{b}\}. 
	\end{align*} 
	This is also a Diophantine set with a description of the form \ref{eqExplFormDioph}. \\
	With \th\ref{DMlaurentseriesmult} the set $d'(\{(a,b,a\cdot b)\mid a,b \in \F_{q^n}\lau{\t}\})$ is a Diophantine set, written in the form of \ref{eqExplFormDioph}. \\
	These are all conditions of \th\ref{H10PeffectiveDioph}, so we get that $d'$ is an effective Diophantine map.\\
	We have that $\F_{q^n}(\t)$ has a recursive representation by \th\ref{DMrecreprKtauDiv}. Furthermore $\F_q(T)$ is a recursive ring, so also $\F_q(T)^n$ has a recursive representation. Since $S_c(\L_{\t,\y_0}\cup\{\a,\b\}) \subset \F_{q^n}\lau{\t} \subset S_c'(\L_{\t,\y_0}\cup\{\a,\b\})$ and $S_c((\L_T\cup\{\b,\TextIn\F_q((T^n))\})^{(n)}) \subset \F_q(T)^n \subset  S_c'((\L_T\cup\{\b,\TextIn\F_q((T^n))\})^{(n)})$, we get that Hilbert's tenth problem exists over \nieuw{$\F_q\lau{T}^n$ and $\F_{q^n}\lau{\t}$}. \nieuw{ \th\ref{DMaddconstFq((T)),DMlaurentseriesAddRel,H10PR2deci+deff=>R1deci}} imply that if Hilbert's tenth problem over \nieuw{$(\F_p(T),\L_T)$} has a positive answer, then Hilbert's tenth problem over $\F_{q^n}\lau{\t}$ has a positive answer.
\end{proof} 

\begin{cor}
	If \th\ref{DMWeakResOfSing}, weak resolutions of singularities in characteristic $p$, holds, then Hilbert's tenth problem over $\F_{q^n}\lau{\t}$ has a positive answer.
\end{cor} 

\begin{proof}
	\nieuw{In \cite{Fq[[T]]ADF} it is proven that Hilbert's tenth problem over $(\F_p\lau{T},\L_T)$ has a positive answer if \th\ref{DMWeakResOfSing} holds. With \th\ref{DMlaurentseriesdecidable} this gives that Hilbert's tenth problem over $\F_{q^n}\lau{\t}$ has a positive answer under the same hypothesis.} 
\end{proof} 

\begin{thm}\th\label{DMlaurentseriesundecidable} 
	We can extend the Carlitz module to a map $\phi: \F_q\lau{T} \ra K\lau{\t}$. Suppose we set on $K$ the special $A$-field structure $\y_0$ of \th\ref{DMdefy0}. Furthermore give $\F_q\lau{T}$ the language $\L_T$ and $K\lau{\t}$ the language $\L_{\t,\y_0}$, then $\phi$ is an effective Diophantine map. \\
	If Hilbert's tenth problem over $(\F_q\lau{T},\L_T)$ has a negative answer, then so has Hilbert's tenth problem over $(K\lau{\t},\L_{\t,\y_0})$. 
\end{thm} 

\begin{proof}
	Since $\phi$ is still injective, we are going to use \th\ref{H10PeffectiveDioph}. \\
	We have that $\t \left(\sum_{i=-N}^\infty a_i\t^i\right) = \sum_{i=-N}^\infty a_i^q \t^{i+1} = \left(\sum_{i=-N}^\infty a_i^q\t^i\right) \t$. Since $a_i^q = a_i$ only holds for $a_i \in \F_q$ we get that
	$$\phi(\F_q\lau{T}) = \F_q\lau{\t} = \{x \in K\lau{\t} \mid x\t = \t x\}.$$
	This is a Diophantine set and we know the explicit formula. \\
	Next we get that $\{\phi(0)\} = \{x\in K\lau{\t}\mid x = 0\}$, $\{\phi(1)\} = \{1\}$ and $\{\phi(T)\} = \{\t\}$ are all Diophantine sets.\\
	We have that $\phi$ is a ring homomorphism, so for $*\in \{+,\cdot\}$ we get that
	\begin{align*} 
		&\phi(\{(a,b,a*b)\mid a,b \in \F_q\lau{T}\}) = \{(\phi(a),\phi(b),\phi(a)*\phi(b))\mid a,b \in \F_q\lau{T}\}\\
		&= \{(x,y,z) \in K\lau{\t}\mid x,y \in \phi(\F_q\lau{T}),\ z = x*y\}. 
	\end{align*} 
	Since $\phi(\F_q\lau{T})$ is a Diophantine set where we know the explicit formulas, the same holds for the sets $\phi(\{(a,b,a*b)\mid a,b \in \F_q\lau{T}\})$ with $*\in \{+,\cdot\}$.\\
	These are all conditions of \th\ref{H10PeffectiveDioph}, so we get that $\phi$ is an effective Diophantine map. \\
	We have that the rings $\F_p(T)$ and $\F_p$ are recursive and by \th\ref{DMrecreprKtau,DMrecreprKtauDiv} this gives that $\F_p(\t)$ is a recursive ring. Since \nieuw{$S_c(\L_T) \subset \F_p(T) \subset S_c'(\L_T)$ and $S_c(\L_{\t,\y_0})\subset \F_p(\t) \subset S_c'(\L_{\t,\y_0})$} we get that Hilbert's tenth problem exists over $(\F_q\lau{T},\L_T)$ and $(K\lau{\t},\L_{\t,\y_0})$. Then with \th\ref{H10PR2deci+deff=>R1deci} and contraposition, we get that if Hilbert's tenth problem over $(\F_q\lau{T},\L_T)$ has a negative answer, then so has Hilbert's tenth problem over $(K\lau{\t},\L_{\t,\y_0})$.
\end{proof} 

\subsection{Equivalence of Models} \label{parDMModels} \phantom{=} \\
In paragraph \ref{parDMKt} we have proven that the Carlitz module $\phi:\F_q[T] \ra K\{\t\}$ is a Diophantine map. If $\phi$ is an arbitrary Drinfeld module with $\rank(\End_K(\phi)) = 1$ it is a Diophantine map. The proof goes the same. This gives different models of $\F_q[T]$ in $K\{\t\}$. In this paragraph, we will give the definition of Diophantine models and their equivalence from \cite{ArtikelDioph} in the case we need. The general case can be found in \cite{ArtikelDioph}. We will then show that the models with different $\phi$ are equivalent if we add the constant $\b$ to the language on $K\{\t\}$. We will do this by constructing Diophantine automorphisms of $K\{\t\}$ and $\F_q[T]$ that take one Drinfeld module to the other.

\begin{defi}
	Let $(R_1,\L_1)$ and $(R_2,\L_2)$ be sets with a language. A \textbf{Diophantine model} of $R_1$ in $R_2$ is an injective Diophantine map $d:R_1 \ra R_2$.
\end{defi}

\begin{defi}
	Let $(R,\L)$ be a set with a language. A \textbf{Diophantine automorphism} of $R$ is a bijective map $f:R \ra R$ such that $f$ and $f^{-1}$ are Diophantine maps. 
\end{defi} 

\begin{defi}
	Let $d_1:R_1 \ra R_2$ and $d_2:R_1 \ra R_2$ be models of $(R_1,\L_1)$ in $(R_2,\L_2)$. We say that $d_1$ and $d_2$ are \textbf{similar} if there exists an Diophantine automorphism $f$ of $R_2$ such that $d_1(R_1) = f(d_2(R_1))$. \\
	The models are \textbf{equivalent} if there exist Diophantine automorphisms $f_1$ of $R_1$ and $f_2$ of $R_2$ such that $d_1 = f_2 \circ d_2 \circ f_1$.
\end{defi} 

\begin{assumption}
	In this paragraph, we assume that the language on $K\{\t\}$ is $\L_{\t,\y}\cup \{\b\}$, where $\b$ is an element such that $\F_p(\b) = \F_q$.
\end{assumption}

\begin{lemma} \th\label{DMequivLemma1}
	The set $\{x \in K\{\t\}\mid x\not=0\}$ is Diophantine.
\end{lemma} 

\begin{proof}
	We have that $\F_q\{\t\} = \{x \in K\{\t\}\mid x\t = \t x\}$, so $\F_q\{\t\}$ is a Diophantine subset of $K\{\t\}$. By \th\ref{DMMRDPFqntau} MRDP holds in $(\F_q\{\t\},\L_{\t,\y}\cup\{\b\})$. The set $\{\t^n \mid n \in \N\}$ is a recursive set, so it is a Diophantine subset of $\F_q\{\t\}$. Since $\F_q\{\t\}$ is a Diophantine subset of $K\{\t\}$ and the languages are the same, this implies that $\{\t^n\mid n\in \N\}$ is a Diophantine subset of $K\{\t\}$.\\
	This implies
	\begin{align*}
		\{x \in K\{\t\}\mid x\not=0\} &= \left\{\sum_{i=k}^n a_i\t^i\mid k \in \N,\ a_i \in K,\ a_k\not=0\right\} \\
		&= \left\{a_k\t^k + \sum_{i=k+1}^n a_i\t^i\mid k\in \N,\ a_i \in K,\ a_k\not=0\right\} \\
		&= \left\{a_k\t^k+f\t^{k+1}\mid f \in K\{\t\},\ a_k \in K\{\t\}^\ti \right\}\\
		&= \left\{a_k\cdot g+f\cdot g \cdot \t \mid g \in \{\t^n\mid n \in \N\}, f\in K\{\t\},\ \ex a \in K\{\t\},\ aa_k=1\right\}.
	\end{align*} 
	Since $\t$ is in the language and $\{\t^n \mid n \in \N\}$ is a Diophantine set, this implies that $\{x \in K\{\t\}\mid x \not= 0\}$ is a Diophantine set.
\end{proof}

\begin{lemma} \th\label{DMequivLemma2}
	Let $\phi,\psi: A \ra K\{\t\}$ be Drinfeld modules that may use a different $A$-field structure on $K$. If $\phi(A) \bs \psi(A)$ is non-empty, it is an infinite set.
\end{lemma} 

\begin{proof}
	We will first prove that $K\{\t\} \cap \psi(F) = \psi(A)$ in $K(\t)$.\\
	Let $x \in K\{\t\} \cap \psi(F)$. Then we can write $x = \frac{\psi(a)}{\psi(b)}$ with $\gcd(a,b)=1$. Since $K\{\t\}$ has a right division algorithm, there exists a $k$ such that $k\psi(b) = \psi(a)$, as otherwise $x \notin K\{\t\}$. Since $\gcd(a,b)=1$, there exist $c,d \in \F_q[T]$ such that $ca+db=1$. Applying $\psi$ gives $\psi(c)\psi(a)+\psi(d)\psi(b)=1$. This means that the two-sided ideal $(\psi(a),\psi(b))$ is $1$. But $\psi(a) =k\psi(b)$, so $1 = (\psi(a),\psi(b)) = (\psi(b))$, which implies that $\psi(b) \in K\{\t\}^\ti = K^\ti$. This gives that $x \in \psi(A)$. Since the other direction follows directly, we get $K\{\t\} \cap \psi(F) = \psi(A)$.\\
	Now take $x \in \phi(A) \bs \psi(A)$, which is possible since it is non-empty. Then for all $n \in \N$, it holds that $x^n \in \phi(A)$. Suppose that $x^{2n}$ and $x^{2n+1}$ are both elements of $\psi(A)$. Then $x = \frac{x^{2n+1}}{x^{2n}} \in \psi(F)$, since $\psi(F)$ is a field. On the other hand $x \in \psi(A) \subset K\{\t\}$, so we get that $x \in \psi(F) \cap K\{\t\} = \psi(A)$. This gives a contradiction, so at least one of $x^{2n}$ and $x^{2n+1}$ is not in $\psi(A)$ and therefore in $\phi(A) \bs \psi(A)$. Since this holds for all $n \in \N$, we get that $\phi(A) \bs \psi(A)$ is infinite.
\end{proof}

\begin{lemma} \th\label{DMequivLemma3}
	Let $\phi$ and $\psi$ be Drinfeld modules, that may use a different $A$-field structure on $K$. Suppose that $\rank(\End_K(\phi)) = \rank(\End_K(\psi)) = 1$. If $\phi(A) \bs \psi(A)$ and $\psi(A) \bs \phi(A)$ are non-empty sets, then there exists a Diophantine automorphism of $K\{\t\}$ such that $d(\phi(A)) = \psi(A)$ and $\phi$ and $\psi$ are similar.
\end{lemma} 

\begin{proof}
	We have that 
	\begin{align*}
		K\{\t\} \bs (\phi(A) \cup \psi(A))  &= (K\{\t\}\bs \phi(A)) \cap (K\{\t\}\bs \psi(A)) \\
		&= \{x \in K\{\t\} \mid x\phi_T - \phi_Tx \not= 0,\ x\psi_T - \psi_Tx \not= 0\}.
	\end{align*}
	With \th\ref{DMequivLemma1} this is a Diophantine set. \\
	The sets $\phi(A) \bs \psi(A) = \{x \in K\{\t\} \mid x\phi_T = \phi_T x,\ x\psi_T-\psi_Tx\not=0\}$ and $\psi(A) \bs \phi(A)$ are also Diophantine by \th\ref{DMequivLemma1}. Then the set $\phi(A) \cap \psi(A)$ is also Diophantine. We will define $d$ separately on these sets.\\ 
	With \th\ref{DMequivLemma2} the sets $\phi(A) \bs \psi(A)$ and $\psi(A) \bs \phi(A)$ are infinite. They are both recursive sets, so let $f$ injectively and recursively enumerate $\phi(A) \bs \psi(A)$ and let $g$ injectively and recursively enumerate $\psi(A) \bs \phi(A)$. Then $f,g,f^{-1}$ and $g^{-1}$ are recursive functions. Since $\rank(\End_K(\phi)) = \rank(\End_K(\psi)) = 1$, we get that the sets $\psi(A)$ and $\phi(A)$ are isomorphic to $A$ as rings. By \cite{FqtMRDP}, MRDP holds in $A$, so also in $\phi(A)$ and $\psi(A)$. This implies that $f,\ g,\ f^{-1}$ and $g^{-1}$ are Diophantine maps.\\
	Now define $d:K\{\t\}\ra K\{\t\}$ by
	$$d(x) = \begin{cases} x &\text{ if }x \in K\{\t\}\bs(\phi(A)\cup \psi(A)) \\ x &\text{ if } x \in \phi(A)\cap \psi(A)\\
		g(f^{-1}(x)) &\text{ if } x \in \phi(A)\bs \psi(A)\\
		f(g^{-1}(x)) &\text{ if }x \in \psi(A)\bs \phi(A).\\ \end{cases} $$
	We have that $K\{\t\}$ is a domain and $\L_{\t,\y}\cup\{\b\}$ is the language of rings with some extra constants. Furthermore $id,\ g\circ f^{-1}$ and $f\circ g^{-1}$ are Diophantine functions and we do a case distinction on Diophantine sets. This means we can use Lemma 5.2.9 from \cite{ArtikelDioph} to get that $d$ is a Diophantine map. By construction we have that $d$ is bijective and $d^{-1} = d$, so $d$ is a Diophantine automorphism.\\ We constructed $d$ in such a way that $d(\phi(A)) = (\phi(A)\cap \psi(A))\cup (\psi(A)\bs \phi(A)) = \psi(A)$, so $\phi$ and $\psi$ are similar.
\end{proof}

\begin{lemma} \th\label{DMequivLemma4}
	Suppose that $\y(T)$ is of finite order. Let $\phi$ and $\psi$ be Drinfeld modules defined over $S_c'(\L_{\t,\y}\cup \{\b\})$. Suppose that $\rank(\End_K(\phi)) = \rank(\End_K(\psi)) = 1$. Then $\phi$ and $\psi$ are of rank $1$ and they are similar.
\end{lemma}

\begin{proof}
	We have that $\phi$ and $\psi$ are defined over $\F_p(\b,\y(T)) = \F_{q^n}$ for some $n\in \Z$. This is a finite field, so with Theorem 4.1.3 of \cite{DM} we get that $\rank(\End_K(\phi)) \geq \rank(\phi)$, so $\rank(\phi)=1$ and, similarly, $\rank(\psi) = 1$. \\
	Suppose that $\phi(A) \subset \psi(A)$. Since $\phi$ and $\psi$ both have rank $1$, their images contain the same amount of elements of each $\t$-degree. This implies that $\psi(A) \subset \phi(A)$. In this case, we can take the identity and they are similar.\\
	So suppose that $\phi(A) \not \subset \psi(A)$. By symmetry, we also get that $\psi(A) \not \subset \phi(A)$. This means that $\phi(A) \bs \psi(A)$ and $\psi(A) \bs \phi(A)$ are both non-empty sets. With \th\ref{DMequivLemma3} we get that $\phi$ and $\psi$ are similar.
\end{proof}

\begin{lemma} \th\label{DMequivLemma5}
	Suppose that $\y(T)$ is of infinite order. Let $\phi,\psi:A \ra K\{\t\}$ be Drinfeld modules with $\rank(\End_K(\phi)) = \rank(\End_K(\psi)) = 1$. Then $\phi$ and $\psi$ are similar. 
\end{lemma}

\begin{proof}
	We start by creating a new $A$-field structure on $K$ with the map $\y_k: \F_q[T] \ra K$ defined by $\y_k(\sum a_iT^i) = \sum \y(a_i)\y(T)^{ki}$. With this we define the Carlitz module $\phi_k$ corresponding to $\y_k$ by $\phi_k(T) = \y_k(T) + \t$. Then $\phi_k(T)$ is defined over $\L_{\t,\y}$ and $\rank(\End_K(\phi_k)) = 1$. \\
	We will prove that $\phi_k(A) \cap \phi_l(A) = \F_q$ if $k\not=l$. Suppose that $x \in \phi_k(A)\cap \phi_l(A)$ and $\deg_\t(x)\geq 1$. Then $x =\phi_k(\sum_{i=0}^n a_i\t^i) = \phi_l(\sum_{i=0}^m b_i\t^i)$. Comparing the degree in $\t$ gives that $n=m$. If we then compare the constant coefficients, we get 
	$$\y_k(\sum_{i=0}^n a_iT^i) = \sum_{i=0}^n \y(a_i)\y(T)^{ki} = \sum_{i=0}^n \y(b_i)\y(T)^{li} = \y_l(\sum_{i=0}^n b_iT^i).$$ 
	Comparing the degree in $\y(T)$ then gives $nk=nl$. We have that $n\geq 1$, since $\deg_{\t}(x)\geq 1$, so $k=l$. This gives a contradiction, so $\phi_k(A) \cap \phi_l(A) \subset \F_q$. Since $\phi_k(c) = \y(c) = \phi_l(c)$ for all $c \in \F_q$, we get $\phi_k(A) \cap \phi_l(A) = \F_q$.\\
	We have that $\phi_k(A)$ contain $q^{n+1}$ elements of $\t$-degree $n$ for all $k \in \Z_{>0}$. Also $\phi(A)$ and $\psi(A)$ both contain at most $q^{n+1}$ elements of $\t$-degree $n$. This implies that there exists a $k \in \Z_{>0}$ such that $\phi(A) \not \subset \phi_k(A)$, $\phi_k(A) \not \subset \phi(A)$, $\psi(A) \not \subset \phi_k(A)$ and $\phi_k(A) \not \subset \psi(A)$. Then with \th\ref{DMequivLemma3} there exist Diophantine automorphisms $d_1,d_2: K\{\t\} \ra K\{\t\}$ such that $d_1(\phi(A)) = \phi_k(A)$ and $d_2(\phi_k(A)) = \psi(A)$. Then $d_2 \circ d_1$ is also a Diophantine automorphism and $d_2(d_1(\phi(A))) = \psi(A)$, so $\phi$ and $\psi$ are similar.  
\end{proof} 

\begin{thm} \th\label{DMequivModels}
	Let $\phi$ and $\psi$ be Drinfeld modules with $\rank(\End_K(\phi)) = \rank(\End_K(\psi)) = 1$. Then $\phi$ and $\psi$ are equivalent models. 
\end{thm} 

\begin{proof}
	We have that $\y(T)$ has finite or infinite order, so with \th\ref{DMequivLemma4} or \th\ref{DMequivLemma5} we get that $\phi$ and $\psi$ are similar. Let $d$ be a Diophantine automorphism of $K\{\t\}$ such that $\phi(A) = d(\psi(A))$. Then $d(\psi(A)) \subset \phi(A)$, so the map $f = \phi^{-1}\circ d\circ \psi$ is well-defined since $\phi$ is injective. By construction of $d$ we can restrict it to a map $d': \F_q(\y(T))\{\t\} \ra \F_q(\y(T))\{\t\}$ and $d'$ is bijective and recursive. Since inverses of recursive functions are recursive, we get that $f$ is a recursive function. Since MRDP holds in $A$, this implies that $f$ is a Diophantine map. By symmetry, we get that $f^{-1}$ exists and is a Diophantine map. This means that $f$ is a Diophantine automorphism. By definition we have $f = \phi^{-1} \circ d \circ \psi$, so $\phi \circ f = d \circ \psi$, which implies $d^{-1} \circ \phi \circ f = \psi$. This gives that $\phi$ and $\psi$ are equivalent.  
\end{proof}

\section{Hilbert's Tenth Problem over Differential Polynomials} \label{chapterDiff}

A variant of the ring of twisted polynomials which also works in characteristic zero is the ring of differential polynomials. The analog of Drinfeld modules in this setting are Krichever modules. In this section, we will prove analogs of some of the results of Section \ref{chapterDM} in this new setting.

\subsection{A Negative Answer over $K[\dd]$} \phantom{=}\\
In this paragraph, we prove that Hilbert's tenth problem over $(K[\dd],\L_R\cup\{\dd\})$ has a negative answer by using the map $\phi:\C[t] \ra K[\dd]$ defined by $\phi(t)= \dd$. Note that this map is a special case of a Krichever module, just as we used the Carlitz module instead of a general Drinfeld module.

\begin{notation}
	We will denote the language $\L_R\cup\{\dd\}$ by $\L_\dd$.
\end{notation} 

\begin{thm}
	Let $(K,\frac{d}{dx})$ be a differential field with field of constants $C$. Give $K[\dd]$ the language $\L_\dd$ and $C[t]$ the language $\L_t$. Define $\phi: C[t] \ra K[\dd]$ by $\phi(\sum_{i=0}^n a_it^i) = \sum_{i=0}^n a_i \dd^i$. Then $\phi$ is an effective Diophantine map and Hilbert's tenth problem over $(K[\dd],\L_\dd)$ has a negative answer.  
\end{thm} 

\begin{proof}
	We have that $\phi$ is an injective map, so we are going to use \th\ref{H10PeffectiveDioph}.\\
	Let $f = \sum_{i=0}^n a_i\dd^i \in K$. Then $f \dd = \sum_{i=0}^n a_i\dd^{i+1}$ and $\dd f = \frac{d}{dx} a_0 + \sum_{i=0}^n (a_{i}+\frac{d}{dx}a_{i+1})\dd^{i+1}$. If $f\dd = \dd f$, this gives that $\frac{d}{dx} a_i = 0$ for all $i$, so $f \in C[\dd]$. This implies
	$$\phi(C[t]) = C[\dd] = \{f \in K[\dd] \mid f\dd = \dd f \},$$
	so $\phi(C[t])$ is a Diophantine set.\\
	We have that $\phi$ is a ring homomorphism, so for $* \in \{+,\cdot\}$ we get
	$$\phi(\{(f,g,f*g)\mid f,g \in C[t]\}) = \{(\phi(f),\phi(g),\phi(f)*\phi(g))\mid f,g \in C[t]\} = \{(f,g,f*g)\mid f,g \in \phi(C[t])\}.$$
	Since we have proven that $\phi(C[t])$ is a Diophantine set, this implies that $\phi(\{(f,g,f*g)\mid f,g \in C[t]\})$ is also a Diophantine set. \\
	Finally we have that $\phi(0) = 0$, $\phi(1)=1$ and $\phi(t) = \dd$, which are all constants of $\L_\dd$. \\
	All sets are given by explicit formulas, so with \th\ref{H10PeffectiveDioph} we get that $\phi$ is an effective Diophantine map.\\
	By \cite{DenefRTchar0} or \cite{DenefRTcharp}, Hilbert's tenth problem over $(C[t],\L_t)$ has a negative answer. We have that $\F_p[t]$ and $\Q[t]$ have a recursive representation. Since $\F_p$ and $\Q$ consist of constants, $\Q[t]$ and $\F_p[t]$ are isomorphic to $\F_p[\dd]$ and $\Q[\dd]$, respectively, so also these rings have a recursive representation. Since we have either $S_c(\L_\dd) \subset \F_p[\dd] \subset S_c'(\L_\dd)$ or $S_c(\L_\dd) \subset \Q[\dd] \subset S_c'(\L_\dd)$, depending on the characteristic, we have that Hilbert's tenth problem exists over $(K[\dd],\L_\dd)$. With \th\ref{H10PR2deci+deff=>R1deci} this implies that Hilbert's tenth problem over $(K[\dd],\L_\dd)$ has a negative answer.
\end{proof} 

\subsection{Results for $K(\dd_1,\ldots,\dd_k)$} \phantom{=}\\
We first prove that $K[\dd_1,\ldots,\dd_k]$ satisfies the Ore condition, so that $K(\dd_1,\ldots,\dd_k)$ exist. Then we show that $\phi: \C(t_1,\ldots,t_k) \ra K(\dd_1,\ldots,\dd_k)$ defined by $\phi(t_i) = \dd_i$ is an effective Diophantine map. We then conclude that Hilbert's tenth problem over $K(\dd_1,\ldots,\dd_k)$ has a negative answer for $k\geq 2$ if the field of constants of $K$ is $\C$ and has a negative answer for $k\geq 1$ if the field of constants of $K$ is $\R$.

\begin{lemma} \th\label{diffOreCondition} 
	Let $K$ be a differential field with derivations $\frac{d}{dx_1},\ldots,\frac{d}{dx_k}$. Then $K[\dd_1,\ldots,\dd_k]$ satisfies the left Ore condition.
\end{lemma} 

\begin{proof}
	Let $a = \sum_{\vec{i}=0}^{\vec{n}} a_{\vec{i}} \dd^{\vec{i}}$ and $b = \sum_{\vec{i}=0}^{\vec{m}} b_{\vec{i}} \dd^{\vec{i}}$ be arbitrary nonzero elements of $K[\dd_1,\ldots,\dd_k]$. Also take $c = \sum_{\vec{i}=0}^{\vec{N}} c_{\vec{i}} \dd^{\vec{i}}$ and $d = \sum_{\vec{i}=0}^{\vec{M}} d_{\vec{i}} \dd^{\vec{i}}$ in $K[\dd_1,\ldots,\dd_k]$. We have that $c\cdot a = d \cdot b$ is a system of linear equations, one for each $\dd^{\vec{i}}$, in the variables $c_{\vec{i}},d_{\vec{i}}$ and coefficients that depend on the coefficients of $a$ and $b$ and their derivatives. Let $n = \max_{j=1}^k(n_j,m_j)$ and take $M_j = N_j = N$. Then the number of equations, which is the number of different $\dd^{\vec{i}}$, is $\prod_{j=1}^k \max(M_i + n_i+1, N_i+m_i+1) \leq (N+n+1)^k$. The number of variables is $\prod_{j=1}^k (M_i+1) + \prod_{j=1}^k (N_i+1) = 2(N+1)^k$. If $k=1$, we take $N=n$ and we get more variables than equations. If $k\geq 2$, we take $N > \frac{n+1-\sqrt[k]{2}}{\sqrt[k]{2}-1}$. This gives $\sqrt[k]{2}(N+1) > N+n+1$, so $2(N+1)^k > (N+n+1)^k$, so also more variables than equations. The system has a solution, namely $c = d= 0$, so it also has a nonzero solution. Since $a$ and $b$ are nonzero and $K[\dd_1,\ldots,\dd_k]$ is a domain, this implies that $c$ and $d$ are both nonzero and $K[\dd_1,\ldots,\dd_k]$ satisfies the left Ore condition.   
\end{proof} 

\begin{assumption}
	For the rest of this paragraph, $K$ will be a differential field with derivations $\frac{d}{dx_1},\ldots,\frac{d}{dx_k}$ and field of constants $C$. We assume that the language on $K(\dd_1,\ldots,\dd_k)$ is $\L_R \cup \{\dd_1,\ldots,\dd_k\}$ and the language on $C(t_1,\ldots,t_k)$ is $\L_R\cup\{t_1,\ldots,t_k\}$. 
\end{assumption} 

\begin{defi}
	We define the $C$-algebra homomorphism $\phi: C(t_1,\ldots,t_k) \ra K(\dd_1,\ldots,\dd_k)$ by $t_i \ra \dd_i$. This is well-defined since the $\dd_i$ commute.
\end{defi} 

Before we can prove that $\phi$ is an effective Diophantine map, we need some lemmas. The proof strategy of this is an extended and modified version of the proof that the Carlitz module extends to an effective Diophantine map from $\F_q(T)$ to $K(\t)$.

\begin{notation}
	For the rest of this paragraph, we set for $1\leq l \leq k$ the field $K_l = \{c \in K\mid \frac{d}{dx_i} c = 0\}$ and the ring $R_{l} = K_l(\dd_1,\ldots,\dd_{l-1},\dd_{l+1},\ldots,\dd_k)[\dd_{l}]$.
\end{notation} 

\begin{lemma} \th\label{diffdivisionRl} 
	For all $x,\kappa \in R_l$ there exist $y,r \in R_l$ such that $x = \kappa y + r$ and $\deg_{\dd_l}(r) < \deg_{\dd_l}(\kappa)$ or $r=0$. 
\end{lemma} 

\begin{proof}
	We will prove this by induction on $\deg_{\dd_l}(x)$. \\
	\textit{Induction basis:} If $\deg_{\dd_l}(x) < \deg_{\dd_l}(\ka)$, the statement holds with $y= 0$ and $r = x$.\\
	\textit{Induction hypothesis:} Suppose that for all $x$ with $\deg_{\dd_l}(x) \leq n$ there exist $y,r \in R_l$ such that $x = \ka y + r$ and $r=0$ or $\deg_{\dd_l}(r)<\deg_{\dd_l}(\ka)$. \\
	\textit{Induction step:} Let $x$ have degree $n+1$. Then we can write $x= a\dd_l^{n+1} + \O(\dd_l^n)$ and $\ka = b \dd_l^m + \O(\dd_l^{m-1})$ with $a,b \in K_l(\dd_1,\ldots,\dd_{l-1},\dd_{l+1},\ldots,\dd_k)$. I use $\O(\dd_l^i)$ here to denote a term of $\dd_l$-degree at most $i$. This gives 
	\begin{align*}
		\ka b^{-1} a \dd_l^{n+1-m} &= b \dd_l^m b^{-1} a \dd_l^{n+1-m} + \O(\dd_l^{m-1}\cdot \dd_l^{n+1-m}) \\
		&= bb^{-1}\dd_l^ma \dd_l^{n+1} + \O(\dd_l^{m-1}\cdot \dd_l^{n+1-m})\\ 
		&= a\dd_l^{n+1} + \O(\dd_l^n).
	\end{align*} 
	This implies that $x-\ka b^{-1} a \dd_l^{n+1-m}$ has degree at most $n$. With the induction hypothesis, there exist $y',r \in R_l$ such that $x-\ka b^{-1} a \dd_l^{n+1-m} = \ka y'+r$ and $r=0$ or $\deg_{\dd_{l}}(r) < \deg_{\dd_l}(\ka)$. This gives $x = \ka (y'+b^{-1}a \dd_l^{n+1-m}) + r$, which proves the induction hypothesis for $n+1$ and therefore the lemma.
\end{proof} 

\begin{notation}
	Let $K$ be a differential field with derivations $\frac{d}{dx_1},\ldots,\frac{d}{dx_k}$. Let $f  = \sum_{\vec{i}} a_{\vec{i}}\dd^{\vec{i}}$. We will abuse notation to write $\frac{d}{dx_j} f$ for $\sum_{\vec{i}} \frac{d}{dx_j}(a_{\vec{i}}) \dd^{\vec{i}}$. 
\end{notation}

\begin{lemma} \th\label{difftechnical} 
	Fix two elements $A,B \in R_{l+1}$. Let $S_{l+1} = \{x \in R_{l+1} \mid x\dd_{l+1} = \dd_{l+1}x\}$ and set $M = \{x \in R_{l+1} \mid Ax\dd_{l+1} = B\dd_{l+1}x\}$. Then there exists a $\kappa \in R_{l+1}$ such that $M = \kappa S_{l+1}$.  
\end{lemma} 

\begin{proof}
	If $M=0$, we can take $\kappa = 0$, so we may assume without loss of generality that $M\not=0$. \\
	Let $\kappa \in M\bs \{0\}$ be an element with the smallest $\d_{l+1}$-degree. Let $x \in S_{l+1}$, then 
	$$A(\kappa x)\dd_{l+1} = (A \kappa \dd_{l+1})x = (B\dd_{l+1}\kappa)x = B\dd_{l+1}(\kappa x).$$
	This proves that $\kappa x \in M$, so $\kappa S_{l+1} \subset M$.\\
	For the converse, take $x\in M$. Then we have $Ax\dd_{l+1} = B\dd_{l+1}x = Bx\dd_{l+1} + B\frac{d}{dx_{l+1}} x$. This implies that $\deg_{\dd_{l+1}}(A) = \deg_{\dd_{l+1}}(B)$ and, since the degree of $B\frac{d}{dx_{l+1}} x$ is lower than that of $Bx\dd_{l+1}$, that $\deg_{\dd_{l+1}}(A-B) < \deg_{\dd_{l+1}}(A)$. \\
	With \th\ref{diffdivisionRl} there exist $y,r \in R_{l+1}$ such that $x = \kappa y + r$ and $r=0$ or $\deg_{\dd_{l+1}}(r) < \deg_{\dd_{l+1}}(\ka)$. Substituting this into $Ax\dd_{l+1} = B\dd_{l+1}x$ gives
	\begin{align*} 
		A\ka y \dd_{l+1} + Ar \dd_{l+1} = B\dd_{l+1} \ka y + B\dd_{l+1}r = B \dd_{l+1} \ka y + Br \dd_{l+1} + B \frac{d}{dx_{l+1}}r. 
	\end{align*} 
	We can rewrite this as
	\begin{align} 
		A\ka y \dd_{l+1} - B \dd_{l+1} \ka y = (A-B) r \dd_{l+1} + B \frac{d}{dx_{l+1}} r. \label{diffEqTechnical} 
	\end{align} 
	The degree of the right-hand side of this equation is at most 
	\begin{align*}
		&\phantom{=}\max(\deg_{\dd_{l+1}}(A-B)+\deg_{\dd_{l+1}}(r)+1,\ \deg_{\dd_{l+1}}(A)+\deg_{\dd_{l+1}}(r)) \\
		&< \max(\deg_{\dd_{l+1}}(A)+\deg_{\dd_{l+1}}(\ka),\  \deg_{\dd_{l+1}}(A)+\deg_{\dd_{l+1}}(\ka)) \\
		&= \deg_{\dd_{l+1}}(A\ka).
	\end{align*}   
	Since $\ka\in M$, we can rewrite the left-hand side of (\ref{diffEqTechnical}) as $A\ka y \dd_{l+1} - A\ka \dd_{l+1} y = A\ka (y\dd_{l+1} - \dd_{l+1} y)$. This has degree strictly lower than that of $A \ka$, since the right-hand side of (\ref{diffEqTechnical}) does so. This implies that $y\dd_{l+1} - \dd_{l+1} y = 0$. This gives that $y \in S_{l+1}$, so $\ka y \in M$. Since $M$ is closed under addition, this implies that $r \in M$. If $r$ is nonzero, we have $\deg_{\d_{l+1}}(r)< \deg_{d_{l+1}}(\ka)$, which leads to a contradiction with the fact that $\ka$ had the lowest degree. This implies that $r=0$, so $x = \ka y \in \ka S_{l+1}$. 
\end{proof} 

\begin{lemma} \th\label{diffOreCondSplit} 
	The ring $R_l$ satisfies the left Ore condition.
\end{lemma} 

\begin{proof}
	Let $a,b \in R_l$ be nonzero elements. Then we can write $a = \sum_{i=0}^n a_i \dd_l^i$ and $b = \sum_{i=0}^m b_i \dd_l^i$ with $a_i,b_i \in K_l(\dd_1,\ldots,\dd_{l-1},\dd_{l+1},\ldots,\dd_k)$. We can inductively clear the denominators of the $a_i$ and $b_i$ since multiplying by an element of $K_l[\dd_1,\ldots,\dd_{l-1},\dd_{l+1},\ldots,\dd_k]$ does not create new denominators nor does it change the $\dd_{l}$-degree of $a$ and $b$. This means that there exists a $c \in K_l[\dd_1,\ldots,\dd_{l-1},\dd_{l+1},\ldots,\dd_k]\bs \{0\}$ such that $ca,cb \in K_l[\dd_1,\ldots,\dd_{l-1},\dd_{l+1},\dd_k][\dd_l] = K_l[\dd_1,\ldots,\dd_k]$. Since $K_l$ is a differential field, we get by \th\ref{diffOreCondition} that $K_l[\dd_1,\ldots,\dd_k]$ satisfies the left Ore condition. This means that there exist $A,B \in K_l[\dd_1,\ldots,\dd_k] \bs\{0\}$ such that $Aca = Bcb$. Since $Ac,Bd \in R_l \bs\{0\}$, this proves that $R_l$ satisfies the left Ore condition.  
\end{proof} 

\begin{lemma} \th\label{diffbeeldphifraction}
	Suppose that the field of constants of $K$ is $C$. Then
	$$C(\dd_1,\ldots,\dd_k) = \{x \in K \mid \va 1\leq j \leq k,\ x\dd_j = \dd_j x\}.$$
\end{lemma} 

\begin{proof}
	We will prove with induction with respect to $l$ that 
	$$\{x \in K(\dd_1,\ldots,\dd_k)\mid \va 1 \leq j \leq l,\ x\dd_j = \dd_j x\} = K_l(\dd_1,\ldots,\dd_k).$$
	\textit{Induction basis:} Since $K_0 = K$, the statement holds for $l=0$.\\
	\textit{Induction hypothesis:} Suppose that the statement holds for some $l \geq 0$.\\
	\textit{Induction step:} With the induction hypothesis we have that
	$$\{x \in K(\dd_1,\ldots,\dd_k) \mid \va 1\leq j\leq l+1,\ x\dd_j = \dd_jx\} 
	= \{x \in K_l(\dd_1,\ldots,\dd_k)\mid x\dd_{l+1} = \dd_{l+1}x\}.$$
	We have with \th\ref{diffOreCondSplit} that $R_{l+1}$ satisfies the left Ore condition. We have that left division rings of fractions are unique. This means we can construct $K_l(\dd_1,\ldots,\dd_k)$ as the left division ring of fractions of $R_{l+1}$. So let $x \in \{y \in K_l(\dd_1,\ldots,\dd_k) \mid y \dd_{l+1} = \dd_{l+1}y\}$. Then we can write $x = \frac{a}{b}$ with $a,b \in R_{l+1}$. With the left Ore condition, there exist $A,B \in R_{l+1}$ such that $Ab\dd_{l+1} = B\dd_{l+1}b$. This implies that $x\dd_{l+1} = \frac{a}{b}\cdot \frac{\dd_{l+1}}{1} = \frac{a\dd_{l+1}}{b}$ and $\dd_{l+1}x = \frac{\dd_{l+1}}{1} \cdot \frac{a}{b} = \frac{B\dd_{l+1}a}{Ab}$. Since $A\cdot b= 1 \cdot Ab$ and $x\dd_{l+1} = \dd_{l+1}x$, we get that $Aa\dd_{l+1} = B\dd_{l+1}a$. Then $a,b \in \{y \in R_{l+1}\mid Ay\dd_{l+1} = B\dd_{l+1}y\} = M$. Define the set $S_{l+1} = \{y \in R_{l+1}\mid x\dd_{l+1} = \dd_{l+1}x\}$. With \th\ref{difftechnical}, there exists a $\ka \in R_l$ such that $M = \ka S_{l+1}$. This means that there exist $a',b'\in S_{l+1}$ such that $x = \frac{\ka a'}{\ka b'} = \frac{a'}{b'}$. We have that \\
	$S_{l+1} = \{y \in R_{l+1} \mid y \dd_{l+1} = \dd_{l+1} y \} = \{y \in R_{l+1} \mid \frac{d}{dx_{l+1}} y = 0\} = K_{l+1}(\dd_1,\ldots,\dd_{l},\dd_{l+2},\ldots,\dd_k)[\dd_{l+1}],$
	which implies that $x \in K_{l+1}(\dd_1,\ldots,\dd_k)$.\\ 
	For the converse, if $x \in K_{l+1}(\dd_1,\ldots,\dd_k)$, then $\frac{d}{dx_j} x = 0$ for all $1 \leq j \leq l+1$. We get $x\dd_j = \dd_j x$ for all $1\leq j \leq l+1$, so $x \in \{x \in K(\dd_1,\ldots,\dd_k)\mid \va 1 \leq j \leq l+1,\ x\dd_j = \dd_j x\}$.\\
	This proves the induction hypothesis for $l+1$, so the statement holds for all $0\leq l \leq k$. \\
	For $l=k$, we get that $K_k$ is the field of constants of $K$, so we have proven the lemma.
\end{proof} 

\begin{thm} \th\label{diffphiFrac}
	The map $\phi: C(t_1,\ldots,t_k) \ra K(\dd_1,\ldots,\dd_k)$ is an effective Diophantine map.
\end{thm} 

\begin{proof}
	We have that $\phi$ is injective, so we will use \th\ref{H10PeffectiveDioph}. \\
	First of all, we have with \th\ref{diffbeeldphifraction} that $\phi(C(t_1,\ldots,t_k))$ is a Diophantine set in the required form.\\
	Secondly we have that $\phi(0) = 0$, $\phi(1) = 1$ and $\phi(t_i) = \dd_i$, which are all constants of the language on $K(\dd_1,\ldots,\dd_k)$. \\
	Finally we have that $\phi$ is a ring homomorphism, so for $* \in \{+,\cdot\}$ we have
	\begin{align*}
		\phi(\{(x,y,x*y)\mid x,y \in C(t_1,\ldots,t_k)\}) &= \{\phi(x),\phi(y),\phi(x)*\phi(y)\mid x,y \in C(t_1,\ldots,t_k)\} \\
		&= \{(x,y,z) \in K(\dd_1,\ldots,\dd_k)\mid x,y \in \phi(C(t_1,\ldots,t_k)),\ z= x*y\}. 
	\end{align*}
	These are Diophantine sets since $\phi(C(t_1,\ldots,t_k))$ is a Diophantine set.\\
	These are all conditions of \th\ref{H10PeffectiveDioph}, so we get that $\phi$ is an effective Diophantine map.
\end{proof} 

To use this theorem in its full generality, we need one more result.

\begin{lemma}[Variant of Lemma 2.2 of \cite{Ct1t2PZSurvey}] \th\label{diffLvsLt} 
	Let $L$ be an infinite field with a language $\L$ that is the language of rings with possibly extra constants. Let $t$ be transcendental over $L$. Then Hilbert's tenth problem over $(L,\L)$ has a positive or negative answer if and only if Hilbert's tenth problem over $(L(t),\L)$ has a positive or negative answer, respectively.   
\end{lemma} 

\begin{cor} \th\label{diffAddti} 
	If Hilbert's tenth problem over $(C(t_1,\ldots,t_l),\L_R\cup\{t_1,\ldots,t_l\})$ has a negative answer, then so has Hilbert's tenth problem over $(C(t_1,\ldots,t_k),\L_R\cup\{t_1,\ldots,t_k\})$ for all $k>l$. 
\end{cor} 

\begin{proof}
	If Hilbert's tenth problem over $(C(t_1,\ldots,t_l),\L_R\cup\{t_1,\ldots,t_l\})$ has a negative answer, we can apply \th\ref{diffLvsLt} $k-l$ times with $L= C(t_1,\ldots,t_{l+i})$ and $\L = \L_R\cup\{t_1,\ldots,t_l\}$ to get that Hilbert's tenth problem over $(C(t_1,\ldots,t_k),\L_R\cup\{t_1,\ldots,t_l\})$ has a negative answer. It then immediately follows that Hilbert's tenth problem over $(C(t_1,\ldots,t_k),\L_R\cup\{t_1,\ldots,t_k\})$ also has a negative answer.
\end{proof} 

\begin{thm}
	Let $K$ be a differential field with field of constants $\C$. Then Hilbert's tenth problem over $K(\dd_1,\ldots,\dd_k)$ with language $\L_R\cup\{\dd_1,\ldots,\dd_k\}$ has a negative answer if $k\geq 2$.
\end{thm} 

\begin{proof}
	First we have that $\Q(\dd_1,\ldots,\dd_k)$ is isomorphic to $\Q(x_1,\ldots,x_k)$, so it has a recursive representation. This implies that Hilbert's tenth problem exists over $K(\dd_1,\ldots,\dd_k)$. \\
	Kim and Roush have proven in \cite{KRC(t1t2)} that Hilbert's tenth problem over $\C(t_1,t_2)$ with language $\L_{t_1,t_2}$ has a negative answer. With \th\ref{diffAddti} the same holds for $\C(t_1,\ldots,t_k)$ with language $\L_R \cup \{t_1,\ldots,t_k\}$. Then with \th\ref{diffphiFrac} and \th\ref{H10PR2deci+deff=>R1deci} we get the result.
\end{proof} 

\begin{thm}
	Let $K$ be a differential field with field of constants $\R$. Then Hilbert's tenth problem over $K(\dd_1,\ldots,\dd_k)$ with language $\L_R \cup \{\dd_1,\ldots,\dd_k\}$ has a negative answer.
\end{thm}

\begin{proof}
	In \cite{DenefRTchar0} Denef proved that Hilbert's tenth problem over $(\R(t),\L_t)$ has a negative answer. With \th\ref{diffAddti} we get that Hilbert's tenth problem over $(\R(t_1,\ldots,t_k),\L_R\cup\{t_1,\ldots,t_k\})$ also has a negative answer. Then by \th\ref{diffphiFrac} and \th\ref{H10PR2deci+deff=>R1deci} Hilbert's tenth problem over $K(\dd_1,\ldots,\dd_k)$ with language $\L_R \cup \{\dd_1,\ldots,\dd_k\}$ has a negative answer.
\end{proof} 

\bibliographystyle{amsplain} 
\bibliography{biblio2}

\end{document}